\documentclass[times,final]{article}
    \makeatletter
    \def\ps@pprintTitle{%
    \let\@oddhead\@empty
    \let\@evenhead\@empty
    \def\@oddfoot{}%
    \let\@evenfoot\@oddfoot}
    \makeatother
    \usepackage{lineno}
    \usepackage{mathrsfs}
    \usepackage[breaklinks=true,pdftex]{hyperref}
    \modulolinenumbers[1]
    \usepackage{amsmath,bm}
    \usepackage{svg}
    \usepackage{comment}
    \usepackage{multirow}
    \usepackage[a4paper,
            bindingoffset=0.1in,
            left=1.5cm,
            right=1.5cm,
            top=1.5cm,
            bottom=1.5cm]{geometry}
    \usepackage{cancel}
    \usepackage{ulem}
    \usepackage{booktabs}
    \usepackage{dirtytalk}
    \usepackage{lineno}
    \usepackage{caption}
    \usepackage{amssymb}
    \usepackage{amsmath}
    \usepackage{savesym}
        \usepackage{algorithmic}
    \usepackage{algorithm}
    \savesymbol{AND}
    \usepackage{adjustbox}
\usepackage{comment}
    \usepackage{rotating}
    \usepackage{graphicx}
    \usepackage{wrapfig}
    \usepackage{amssymb}
    \usepackage{soul,color}
    \usepackage{subcaption}
    \usepackage{lineno}
\usepackage{xcolor}
\usepackage{url}
    \soulregister\cite7
    \soulregister\ref7
    \soulregister\pageref7

    \usepackage{todonotes}

    \DeclareRobustCommand{\uvec}[1]{{%
    \ifcsname uvec#1\endcsname
    \csname uvec#1\endcsname
    \else
    \bm{\hat{\mathbf{#1}}}%
    \fi}}
\mathchardef\breakingcomma\mathcode`\,
{\catcode`,=\active
  \gdef,{\breakingcomma\discretionary{}{}{}}
}

\begin{document}
\begin{center}
{\LARGE \textbf{Accelerating Dimensionality Reduction in Wave-Resistance Problems through Geometric Operators}}\\
\vspace{0.5cm}
{\small Stamatios Stamatatelopoulos$^{1,2,*}$\let\thefootnote\relax\footnote{$^*$Corresponding author. E-mail address: stamatis@mit.edu (S. Stamatatelopoulos)}},
{\small Shahroz Khan$^{1,3}$}
{\small Panagiotis Kaklis$^{1,4}$}
\\\vspace{0.2cm}
{\small $^{1}$Department of Naval Architecture, Ocean and Marine Engineering, University of Strathclyde, Glasgow (UK)}\\
{\small $^{2}$ Department of Mechanical Engineering, Massachusetts Institute of Technology,
Cambridge, (USA)}\\
{\small $^{3}$BAR Technologies Ltd, Portsmouth (UK)}\\
{\small $^4$Foundation for Research \& Technology Hellas (FORTH), Institute of Applied \& Computational Mathematics (IACM), \\ Division: Numerical Analysis \& Computational Science, Group: Data Science, Heraklion, Crete (GR)}\\
\end{center}

\section*{\centering Abstract}
Reducing the dimensionality and uncertainty of design spaces is a key prerequisite for shape optimisation in computationally intensive fluid problems. However, running these analyses at an offline stage itself poses a computationally demanding task. In this work, we propose a unique framework for the inexpensive implementation of sensitivity analyses for reducing the dimensionality of the design space in wave-resistance problems. At the heart of our approach is the formulation of a geometric operator that leverages, via high-order geometric moments, the underlying connection between geometry and physics, specifically the wave-resistance coefficient ($C_w$), of ships using the slender body theory based on the well-known Vossers' integral. The resulting geometric operator is computationally inexpensive yet physics-informed and can act as a geometry-based surrogate to drive parametric sensitivities. To analytically demonstrate the capability of the proposed approach, we use a well-known benchmark geometry, namely, the modified Wigley hull. Its simple analytical formulation allows for closed expressions of the geometric operators and exploration of computational domains that would otherwise be inaccessible. In this context, the proposed geometric operator outperforms existing similar approaches by achieving 100\% similarity with $C_w$ at a fraction of the computational cost.
\vspace{0.2cm}\\
\textit{Keywords:} Dimensionality Reduction; Sensitivity Analysis; Shape Optimisation; Wave Resistance; Geometric Operators

\section{Introduction}

\par\noindent Identifying satisfactory solutions to physics-based optimisation problems, such as wave resistance minimisation, is a task of fundamental importance in the community of Engineering Design. However, the relevant optimisation algorithms may require evaluating the physics-based objective a number of times, which can be computationally prohibitive. In fact, modern design tools can sometimes prove to be too expensive to run even in cases where a relatively small number of responses need to be evaluated \cite{sun2020surrogate}. In tandem with the prohibitive cost of physics-based solvers is the well-known curse of dimensionality \cite{koppen2000curse,khan2018generative,chen2015measuring} which, in the context of sample-based design optimisation, is realised by an explosion in the number of sample evaluations required by optimisation algorithms as the number of optimisation parameters increases. This is precisely the point of intervention of the so-called method of Dimensionality Reduction (DR), whose preliminary  objective is to reduce the dimensionality of the design space while keeping the same design variability as the original design space. 

\par\noindent In the DR-pertinent literature, there are well-studied unsupervised techniques (e.g., Principal Component Analysis (PCA)\cite{ssdr_r32}, auto-encoders \cite{ssdr_r25}, etc.) or supervised approaches (e.g., Active Subspace Method (ASM) \cite{ConstantineActiveSubspaces}, Sensitivity Analysis (SA) \cite{KHANintrasensitivity}), which may or may not require evaluation of designs' physics, respectively. Among these techniques, PCA, auto-encoders, and ASM extract the latent features of the original design space to create a lower-dimensional subspace while capturing the maximum geometric variability. In contrast, SA is a selection process that identifies parameters that are less sensitive (or insensitive) to physics. As parameters with lower sensitivity have a negligible effect on performance, they can be excluded to reduce the dimensionality of the design space \cite{saltelli2008global,saltelli2010avoid,cacuci2003sensitivity}.\\

\par\noindent Compared to SA, unsupervised DR techniques can be computationally efficient since they do not rely on performance labels. However, their effectiveness may be limited when there is no direct correlation with shape modifications. In contrast, SA implementation is more informed because it not only reduces dimensionality but also provides valuable insights into the driving features of designs that contribute to the extreme variability in performance \cite{KHANintrasensitivity,GMGSA}. Consequently, SA enables users to allocate resources more effectively from the early stages of design, thus expediting product development. However, implementing SA can be computationally demanding, especially when an analytical solution is not available, and costly numerical simulations become necessary. Although meta/surrogate models can accelerate SA, constructing the surrogate itself can be a computationally intensive task due to the requirement of evaluating performance labels for dataset creation \cite{cheng2020surrogate}.\\

\par\noindent Therefore, the aim of this study is to reduce the computational cost of SA for DR purposes by leveraging the underlying connection between physics and geometry. This is achieved by considering that evaluating geometry-based operators is much cheaper compared to their physics-based counterparts. More precisely, for the wave-resistance operator, two geometry-based operators are proposed and tested as surrogate models for SA. It is important to note that the only requirement for the surrogate geometric operators is to be sensitive, in a sense that will be defined below, to the same parameters as the physics-based quantity they are related to. This requirement is much less strenuous when compared to surrogate modelling for optimisation, where the surrogate model must mimic the complete behaviour of the high-fidelity physical model.\\ 

\par\noindent To build a vigorous analytical foundation of the proposed approach we used a wigley hull as a test case, which is a standard benchmark geometry in the Naval Architecture domain. It is a simplified mathematical model of ship hulls and is described by simple polynomial representations, making it analytically tractable. The simple analytical expression of the Wigley hull allows for closed forms of geometric operators, thereby reducing computational costs and providing access to computational domains that would otherwise be inaccessible. Leveraging these closed-form solutions at the preliminary design stage, it offers the potential for explicit analysis of the relevant operators, facilitating a deeper understanding of their behaviour and implications.\\

\par\noindent The remainder of this paper is organised as follows: Section \ref{sec:2} provides a brief overview of the concepts of SA and Uncertainty Analysis (UA) which are fundamental to DR. In Section \ref{sec:3}, a general framework is developed for matching physics-based operators to geometry-based ones, two such candidate operators are introduced and the DR framework employed in this study is discussed in more detail. In Section \ref{sec:4}, the parametric modeller and the physics-based property, that the methodology developed in the previous section will be applied to, are established. In Section \ref{sec:5} the relevant results are displayed and analysed. Finally in Section \ref{sec:6}, conclusive remarks and potential next steps are provided.

\section{Background}\label{sec:2}

\par\noindent The design of mathematical models which describe man-made systems reliably is of fundamental importance to the scientific community. The principal role of SA and UA is to answer questions that arise in the design and validation of such models. Such questions can be expressed as (see, e.g., \cite[p.\ 3]{cacuci2003sensitivity}, \cite[p.\ 183]{saltelli2008global}):
\begin{itemize}
    \item \textit{How well does the model under consideration represent the underlying physical phenomena?}
    \item \textit{How far can the calculated results be extrapolated and how can this be achieved? }
    \item \textit{Which factor or group of factors is most responsible for producing model outputs within or outside specified bounds?}
\end{itemize}
Specifically, the role of UA given certain model input conditions (such as probability distributions for the input variables) is to \textit{quantify the resulting effect on the output} of the model. On the other hand, the role of SA is to \textit{identify which input variables contributed most to said effect}. Grounded on the above, a reasonable definition of SA and UA can be stated as:

\begin{quotation}
\par\noindent {\it Sensitivity and Uncertainty Analysis is a methodology for the formal evaluation of data and models \cite[p.\ 3]{cacuci2003sensitivity}. Uncertainty Analysis is the quantification of uncertainty in the model output given uncertainties in the model input. Sensitivity Analysis is the apportioning of the uncertainty in the model output to different sources of uncertainty in the model input \cite[p. 14]{saltelli2008global}.}
\end{quotation}
Even though the goals of these two analyses differ, in practice, they are often performed simultaneously and the term \textit{Sensitivity Analysis} has prevailed \cite[p.\ 1]{saltelli2010avoid}. For the sake of succinctness, this terminology has also been adopted in this paper to refer both to SA and UA.\\

\par\noindent Having provided a definition of SA that is agreed upon in the literature, we proceed by reviewing the various approaches to perform SA. Out of the many dichotomies of such methods\footnote{Other than \textit{Local} vs \textit{Global}, other distinctions include \textit{Statistical} vs \textit{Deterministic} (see \cite[p.\ 8]{cacuci2003sensitivity} for a brief discussion) or \textit{Qualitative} vs \textit{Quantitative} (see \cite{Iooss2015}). In this context, the SA methodology presented in Section \ref{sec:SobolsSA} can be characterised as global, quantitative and deterministic.}, this brief review will focus on the distinction of \textit{Local} versus \textit{Global} SA since the latter is central to this study. Local SA (LSA) investigates the model behaviour in a small neighbourhood of predetermined points of high interest. On the other hand, Global SA (GSA) investigates the model behaviour over the entire parameter domain. \\

\par\noindent Historically, LSA is the predecessor of GSA with the first systematic methodology being formulated in \cite{bode1945network} as reported in \cite[p.\ 4]{cacuci2003sensitivity}. Perhaps the simplest of such methods is the {\it One at A Time} (OAT), where after a \textit{base point} is chosen, each of the parameters is varied by a specified amount, while all other others are held constant. Comparing the resulting change in model output between parameters can provide insight in their relative importance. Symbolically, given a model $f:\mathbb{R}^n \to \mathbb{R}$ one assigns a base point $\mathbb{R}^n\ni{\bf X}^0 = (x^0_1,...,x^0_n)$ and a fixed variation $\epsilon_i\in\mathbb{R}$, $i=1,...,n$ for each parameter. Then, the respective \textit{sensitivity index} for the $i^{th}$ parameter will be $|f({\bf X}_0) - f(x^0_1,...,x^0_i+\epsilon_i,...,x^0_n)|$. One of the shortcomings of OAT is that by fixing all parameters other than the $i^{th}$ one, it is impossible to investigate the interactions between two different parameters. For example take $n=2$: it might be the case that $|f(x^0_1+\epsilon_1,x_2^0+\epsilon_2)-f(x^0_1,x^0_2)|$ is much greater than either of the two sensitivity indexes but it will never be evaluated. The more sophisticated {\it Elementary Effects} method (EE), also known as \textit{Morris} method or \textit{winding stairs} \cite{saltelli2010avoid}, mitigates this issue by performing a number of \textit{steps} where in each step a new OAT experiment is performed \cite{Iooss2015} and thereafter conclusions for each parameter can be drawn by considering all relevant indices.\\

\par\noindent The OAT and EE methods belong in the category of \textit{screening} methods, characterised by the subdivision of the domain of each parameter into a number of levels. For a more detailed review of LSA methods, which includes derivative-based methods, the reader is referred to \cite{BORGONOVO2016SAreview}. As stated earlier, in order to investigate the entire parameter domain one employs GSA methods. Such methods include regression-based, variation-based and density-based approaches among others (see \cite{BORGONOVO2016SAreview}). Frameworks for determining which method to employ have been proposed in the literature \cite{Iooss2015}, with the first step being the determination of how close to linear is the relationship of the model output to the model input. More specifically, for any constants $\beta_i\in\mathbb{R}$, the question is how well does the linear relationship $f({\bf X}) =$ $\beta_1x_1 + ... + \beta_nx_n$
represent the given model. To handle this issue, one can perform a number of tests such as evaluating the so-called \textit{Pearson correlation coefficient} for each parameter to determine whether its relationship with the output is linear. For any parameters that pass the chosen linearity test satisfactorily, there are readily available measures to evaluate the  parameter's sensitivity such as the Standard Regression Coefficient.\\

\par\noindent If the user chooses not to make any assumptions about the model or the linearity hypothesis fails by the procedure above, \cite{Iooss2015} suggest to employ Sobol's variance-based method. Introduced by Ilya M. Sobol \cite{sobol1990sensitivity}, the method is based on the computation of the so called \textit{Sobol's Index} for each of the input parameters. This method applies to scalar models ($f({\bf X})\in\mathbb{R}$), however, there do exist generalisations to more than one outputs (see \cite{GAMBOAMultiOutputs}). Since this is the SA approach that has been adopted for this study, more details are provided in Section \ref{sec:SobolsSA}.\\

\par\noindent We proceed to the concept of DR which relates to SA via the last of the questions in the start of this section; \textit{Which factor or group of factors is most influential to the model output?} Having calculated Sobol's indices $SI_i\in\mathbb{R}$ for all parameter values, a straightforward approach to answering this question is via specifying a threshold value $\epsilon$ \cite{SARRAZIN2016GSA} and filtering the sensitive parameters accordingly. This approach is capable of identifying \textit{subsets} of parameters but fails to determine \textit{linear combinations} of the parameters, which is a feature present in other approaches to DR such as PCA \cite{sorzano2014survey} or ASM \cite{ConstantineActiveSubspaces}.\\

\par\noindent Fundamental to the current study is the work of \cite{GMGSA}. The authors perform SA on the wave resistance coefficient $C_w$ of a series of ship hulls using the composition of a parametric modeler and a $C_w$-solver as the model $f(\cdot)$ under investigation. SA is performed globally via the variance-based method of Sobol's indices as discussed above. Motivated by the well studied correlation between the so-called Sectional Area Curve (SAC)\footnote{SAC represents the longitudinal variation of the cross-sectional area of a ship's hull below the waterline. The integrated area under the SAC gives the displaced volume and provides an effective and simple description of global geometric properties. At the same time, it is closely related to the resistance and propulsion performance of a ship.} of a ship's hull and ship wave resistance (see \cite{han2012hydrodynamic}) the authors construct a geometry-based operator on ship hulls dubbed \textit{Shape Signature Vector} (SSV) and perform an alternative SA by composing the parametric modeler with the SSV. The two SA approaches are studied for correlation with the expectation that the much cheaper to compute SSV would produce results related to the physics-based solver and thus could be used to expedite the preliminary design stage of hull design. The outcome of this study does show a positive correlation between the two SA's and this is the starting point for the work relevant to this study.

\section{Sensitivity Analysis using Geometric Operators}\label{sec:3}

\par
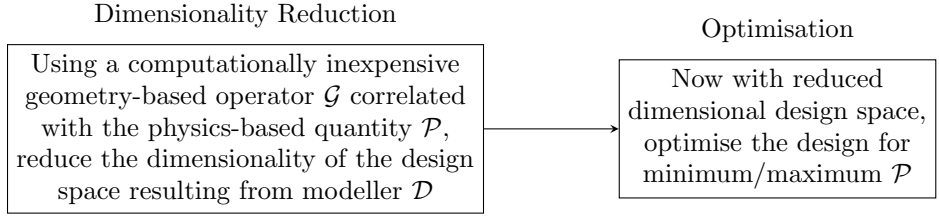
\begin{figure}[H]
    \centering
    \begin{tikzpicture}
        \node[align=center] at (0, 1.5) {Dimensionality Reduction};
        \node[align=center] at (7, 1.3) {Optimisation};

        \node[draw, align=center] (box1) at (0,0) {Using a computationally inexpensive\\
            geometry-based operator ${\cal G}$ correlated\\
            with the physics-based quantity ${\cal P}$,\\
            reduce the dimensionality of the design\\
            space resulting from modeller ${\cal D}$};

        \node[draw, align=center] (box2) at (7,0) {Now with reduced \\
             dimensional design space,\\
            optimise the design for \\minimum/maximum ${\cal P}$};
            
       \draw [-stealth](box1.east) -- (box2.west);
    \end{tikzpicture}
    \caption{Illustration of the application of the geometric operator ${\cal G}$. It's important to note that ${\cal G}$ is not employed in the optimisation step; rather, its use is confined to identifying the most sensitive parameters for dimension reduction. This enables an expedited physics-based optimisation process.}
    \label{fig:DimensionalityReductionoptimisationDiagram}
\end{figure}

\noindent Let be given a parametric modeller ${\cal D}$ of 3D objects embedded in an $n$-dimensional design space ${\cal X} = [a_1,b_1]\times...\times[a_n,b_n]$ where each design parameter $t_i,\ i=1,...,n$ varies continuously in $[a_i,b_i]$. Further, if ${\cal P}({\cal D}) \in \mathbb{R}$ is a physics-based property of objects generated by ${\cal D}$, perform SA on ${\cal D}$ with ${\cal P}$ being the Quantity of Interest (QoI), to identify a subset of $k\le n$ parameters $I=\{i_1,...,i_k\}$ (i.e. $\{t_{i_1},...,t_{i_k}\}$) to which ${\cal P}$ is most sensitive to. We aim at identifying a geometry-based operator ${\cal G}$ such that performing SA with ${\cal G}$ as the QoI, results in a subset of $k'<n$ sensitive parameters $J=\{j_1,...j_{k'}\}$ such that the now reduced-dimensionality design space comprising of 
the parameters in $J$ can be used for optimisation with $\cal P$ with acceptable loss of information. This is illustrated in Figure \ref{fig:DimensionalityReductionoptimisationDiagram}, 

\par
\begin{figure}[H]
    \centering
    \includegraphics[width=0.8\linewidth]{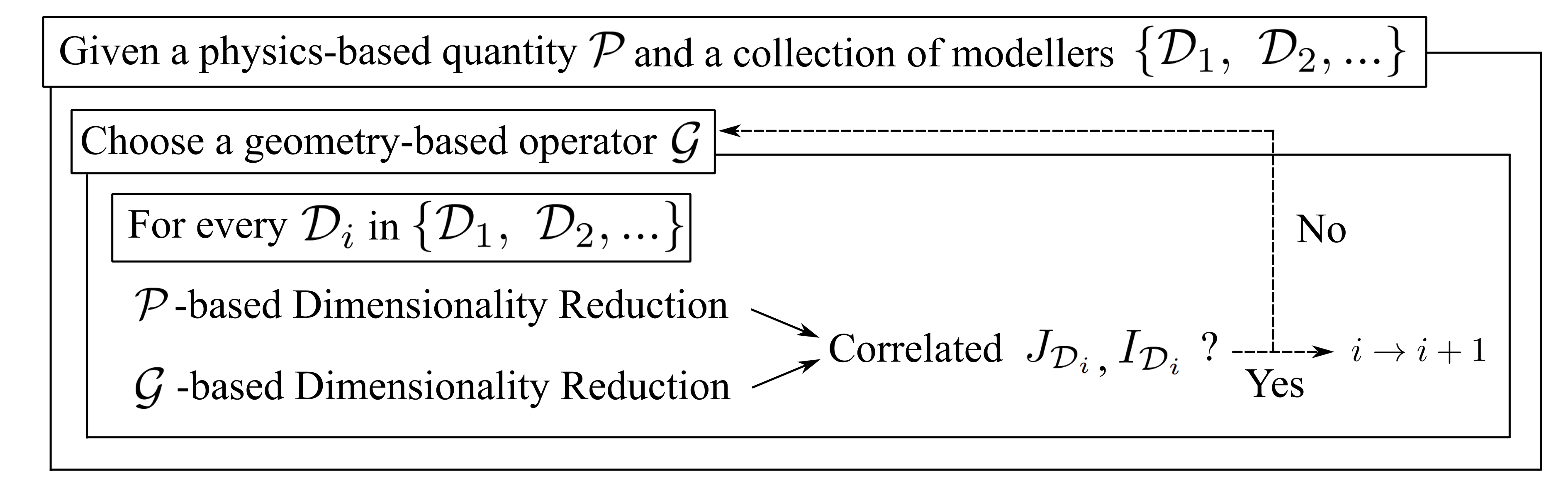}
    \caption{Diagrammatic representation of the methodology for identifying a suitable geometric operator ${\cal G}$ compatible with a physics-based quantity ${\cal P}$ across a set of parametric modellers ${{\cal D}_1, {\cal D}_2, \ldots}$.}
    \label{fig:MethodologyDiagram}
\end{figure}

\noindent It is important to emphasise that this process produces triplets $({\cal D},\ {\cal P},\ {\cal G})$ without any guarantee that the correlation between $\cal P$ and $\cal G$ is invariant to $\cal D$. It is essential therefore to investigate the dependence of the pair $({\cal P},{\cal G})$ on the parametric modeller ${\cal D}$ for, if this dependence is weak, ${\cal G}$ can readily be applied to a range of similar parametric modellers. In this effort, the process above can be repeated for constant $\cal P$ and $\cal G$ across a set of various parametric modellers $\{{\cal D}_1,{\cal D}_2,...\}$. By verifying the compatibility of an increasing amount of triplets $({\cal D}_i,\ {\cal P},\ {\cal G})$, the dependence of $({\cal P},\ {\cal G})$ on the parametric modeller is deceased. For example, in the present case, the aim is to identify a geometric operator ${\cal G}$ which is compatible with ${\cal P}$ equal to the wave-making resistance coefficient $C_w$ for a set of parametric modelers $\{{\cal D}_1,{\cal D}_2,...\} = \{{\cal D}_i:\ \text{${\cal D}_i$ produces slender hulls}\}$. The above methodology for identifying triplets $(\{{\cal D}_1,{\cal D}_2,...\},{\cal P},{\cal G})$ is outlined in Figure \ref{fig:MethodologyDiagram}.

\par\noindent In this paper, the procedure outlined in Fig.\,\ref{fig:MethodologyDiagram} has been applied twice, for two different geometric operators, on a single parametric modeller ${\cal D}_1$. Both of these geometric operators are based on the concept of geometric moments, introduced in the following three sections. Then, the chosen approach for DR as well as the approach for performing correlation analysis between $J_{{\cal D}_i}$ and $I_{{\cal D}_i}$ are outlined.

\subsection{Geometric Moments}

\par\noindent There is a plethora of geometric operators $\cal G$ which have been studied in the literature, such as the so-called Fourier Descriptors (FD) \cite{zhang2002comparative}, which extract information from the geometry of the contour, Convexity Measures (CM) (\cite{rahtu2006new};  \cite{corcoran2011convexity}) or even Elongation Measures \cite{stojmenovic2008measuring}, which at their simplest form are defined as the width over the length of the bounding box of a 2D shape. In order to apply the procedure in Fig.\,\ref{fig:MethodologyDiagram}, two desirable attributes on $\cal G$ are the following:
\begin{subequations}
\begin{gather}
    \text{$\cal G$ is able to extract enough information from the underlying geometry;}\label{eq:Grequirements-1}\\
    \text{$\cal G$ is easy and fast to evaluate for arbitrary geometries.}\label{eq:Grequirements-2}
\end{gather}
\end{subequations}
In view of the above two points, the two operators $\cal G$ which have been investigated in this study are based on geometric-moments $M(p,q,r)$. Starting with the definition of $M(p,q,r)$, the $s^{th}$ order moment where $s=p+q+r$ over the solid ${\cal D}$ is defined as
\begin{gather}\label{eq:Mpqr}
    M(p,q,r) = \iiint_{\cal D}x^py^qz^rdxdydz.
\end{gather}
From the study of the renowned Problem of Moments \cite{shohat1950problem} it is well known that the geometric moments (\ref{eq:Mpqr}) contain significant information regarding the underlying geometry, which can be used to approximate it and, under certain conditions, even uniquely reconstruct it \cite{cuyt2005multidimensionalintegralinversion,milanfar2000shape,kousholt2021reconstruction,ghorbel2005reconstructing}. For example, the volume of $\cal D$ is given by $V_{\cal D}=M(0,0,0)$ and its centroid $C = (C_x,C_y,C_z)$ by,
\begin{subequations}\label{eq:centroid}
\begin{align}
    C_x &= \frac{M(1,0,0)}{M(0,0,0)},\quad C_y=\frac{M(0,1,0)}{M(0,0,0)},\quad C_z=\frac{M(0,0,1)}{M(0,0,0)}.
\end{align}
\end{subequations}
Further, there are numerous procedures for evaluating moments of arbitrary shapes efficiently, one of which is presented in \cite{pozo2010efficient} and has been \href{https://github.com/Stamatis8/geom_moments}{implemented} in the scope of this study. It is then evident that geometric-moments easily satisfy (\ref{eq:Grequirements-1}) and (\ref{eq:Grequirements-2})\\

\par\noindent The utility of moments as a mathematical construction has been recognised in various disciplines, with applications starting from probability and statistics \cite{diaconis1987application} and including image analysis \cite{teh1988image}, signal processing \cite{sezan1987incorporation}, feature extraction \cite{khan2022shape}, computing tomography \cite{milanfar1995reconstructing,milanfar1996moment} and inverse potential theory \cite{strakhov1986uniqueness,brodsky1990concerning} among others. Most importantly, however, they have shown promising results in a context  similar to this study \cite{GMGSA,khan2023shiphullgan} via construction of so called Shape Signature Vector (SSV).

\subsection{Shape Signature Vector}\label{sec:SSV}
\par\noindent As mentioned above, significant geometric information carried by ${\cal D}$ can be captured with a large enough number of its moments. It is then straightforward to consider as a potential moment-based geometric operator, the vector consisting of all moments up to a certain order. This geometric-based operator is called Shape Signature Vector (SSV) and was first introduced in \cite{GMGSA} where it was investigated in an analogous context to the present study. Before introducing the relevant notation, notice that the moments in (\ref{eq:Mpqr}) depend on the position of $\cal D$ in its ambient space and its scale, two non-desirable characteristics. A translation-invariant form of $M(p,q,r)$ can readily be defined if ${\cal D}$ is positioned at its centroid,
\begin{gather}\label{eq:MTpqr}
    M_T(p,q,r) = \iiint_{\cal D}(x-C_x)^p(y-C_y)^q(z-C_z)^rdxdydz.
\end{gather}
However, (\ref{eq:MTpqr}) still varies when ${\cal D}$ is scaled uniformly by $\lambda$, for if we apply the transformation $g(x,y,z)=(\lambda x, \lambda y, \lambda z)$ to (\ref{eq:Mpqr}) over the scaled geometry $\lambda {\cal D}$,
\begin{align}
    \iiint_{\lambda{\cal D}}x^py^qz^rdxdydz &= \iiint_{{\cal D}}(\lambda x)^p(\lambda y)^q(\lambda z)^r|\text{det}D_g|dxdydz\nonumber\\
    &= \iiint_{{\cal D}}(\lambda x)^p(\lambda y)^q(\lambda z)^r\lambda^3dxdydz\nonumber\\
    &=\lambda^{p+q+r+3}M(p,q,r)\label{eq:MpqrScaled}.
\end{align}
Then, looking at (\ref{eq:centroid}) and (\ref{eq:MpqrScaled}) we can write,
\begin{align}
    \iiint_{\lambda{\cal D}}(x-C_x)^p(y-C_y)^q(z-C_z)^rdxdydz &= \iiint_{\cal D}(\lambda x-\lambda C_x)^p(\lambda y-\lambda C_y)^q(\lambda z-\lambda C_z)^r\lambda ^3dxdydz\nonumber\\
    &= \lambda^{p+q+r+3}M_T(p,q,r).
\end{align}
Since for $\lambda {\cal D}$, $M(0,0,0)=M_T(0,0,0)=V_{\lambda {\cal D}}=\lambda^3V_{\cal D}$, a scaling invariant moment $M_S(p,q,r)$ and the translation \& scaling invariant moment $M_I(p,q,r)$ of ${\cal D}$ are given by,
\begin{align}
    M_S(p,q,r) &= \frac{M(p,q,r)}{M(0,0,0)^{(p+q+r+1)/3}}\label{eq:MSpqr},\\
    M_I(p,q,r) &= \frac{M_T(p,q,r)}{M_T(0,0,0)^{(p+q+r+1)/3}}\label{eq:MIpqr}.
\end{align}
Now, given a solid ${\cal D}$, consider the set $M_I^s$ of all scaling and translation invariant moments of ${\cal D}$, $M_I(p,q,r)$ of order $s$.
\begin{gather}
    M_I^s = \{M_I(p,q,r):\ p+q+r = s\}
\end{gather}
For example, 
\begin{align*}
    M_I^0 &= \{M_I(0,0,0)\}\\
    M_I^1 &= \{M_I(1,0,0),\ M_I(0,1,0),\ M_I(0,0,1)\}\\
    M_I^2 &= \{M_I(2,0,0),\ M_I(1,1,0),\ M_I(1,0,1),\ M_I(0,2,0),\ M_I(0,1,1),\ M_I(0,0,2)\}.
\end{align*}
Notice that all translation invariant moments of order 1 are equal to zero by definition. In light of this, define the Shape Signature Vector (SSV) of order $N$, as the vector consisting of all elements of $M_I^s$ for $s=0,2,...,N$,
\begin{gather}\label{eq:SSV}
    SSV_N = [M_I^0,M_I^2,...,M_I^N].
\end{gather}
We conclude this section by calculating the cardinality of $SSV_N$. For each set $M_I^s$, Card($M_I^s$) = $(s+1)(s+2)/2$ which can be proven inductively. Then since all $M_I^s$ are disjoint,
\begin{align}\label{eq:SSVcardinality}
    \text{Card}(SSV_N)=\sum_{s=0,\ s\ne 1}^N\text{Card}(M_I^s)= \sum_{s=0,\ s\ne 1}^N\frac{(s+1)(s+2)}{2}= \frac{1}{6}N^3+N^2+\frac{11}{6}N-2,\ N\ge 1.
\end{align}

\subsection{Slender Body Operator}\label{sec:SlenderBodyOperator}

\par\noindent As will be explained in Section \ref{sec:4}, the chosen physics-based quantity in this study is the wave resistance coefficient $C_w$. The well known correlation between $C_w$ and the SAC of a hull is the starting point for the derivation of the slender body operator. More specifically, in the context of slender body theory, Vossers proposed the following expression for estimating the wave resistance of a slender ship moving at constant velocity $U$, at water density $\rho$, with $K=g/U^2$, for $g$ gravitational acceleration  \cite{vossers1962some},
\begin{gather}\label{eq:vossers-integral-full}
    \frac{R}{-0.5\rho U^2}=\iint_{\mathbb{R}^2} S''(x)S''(\xi)Y_0(K|x-\xi|)dS'(x)dS'(\xi)
\end{gather}
where $Y_0$ is the Bessel function of the second kind, $S(x)$ is the SAC of the ship's hull with $S(x)=0$ for $x\notin[-L/2,L/2]$ and integration is of Riemann–Stieltjes type. Next, under the assumption that $S'(\pm L/2) = 0$, (\ref{eq:vossers-integral-full}) transforms to the following Riemann integral, 
\begin{gather}\label{eq:vossers-integral}
    \frac{R}{-0.5\rho U^2}=\int_{-L/2}^{L/2}\int_{-L/2}^{L/2}S''(x)S''(\xi)Y_0(K|x-\xi|)dxd\xi,
\end{gather}
Recall the definition of $Y_0(x)$ to be
\begin{gather}
    Y_0(x) = \frac{2}{\pi}(\ln(0.5x)+\gamma)J_0(x)+\frac{2}{\pi}\sum_{k=1}^\infty(-1)^{k-1}\frac{(0.25x^2)^k}{(k!)^2}\sum_{j=1}^k\frac{1}{j}\label{eq:Y0},\\
    J_0(x) = \sum_{k=0}^\infty\frac{(-0.25x^2)^k}{(k!)^2},
\end{gather}
with $J_0$ the Bessel function of the first kind. We can then write,
\begin{align}
    Y_0(K|x-\xi|) &= \frac{2}{\pi}(\ln(0.5K|x-\xi|)+\gamma)J_0(K|x-\xi|)+\frac{2}{\pi}\sum_{k=1}^\infty(-1)^{k-1}\frac{(0.25(K|x-\xi|)^2)^k}{(k!)^2}\sum_{j=1}^k\frac{1}{j}\nonumber\\
    &= \frac{2}{\pi}\ln(|x-\xi|)J_0(K|x-\xi|) + \frac{2}{\pi}(\ln(0.5K)+\gamma)J_0(K|x-\xi|)\nonumber\\
    &\hspace{10em}+\frac{2}{\pi}\sum_{k=1}^\infty(-1)^{k-1}\frac{K^{2k}}{2^{2k}(k!)^2}\Bigg(\sum_{j=1}^k\frac{1}{j}\Bigg)(x-\xi)^{2k}\nonumber\\
    &= \frac{2}{\pi}\ln(|x-\xi|)J_0(K|x-\xi|) + \frac{2}{\pi}(\ln(0.5K)+\gamma)\sum_{k=0}^\infty(-1)^{k}\frac{K^{2k}}{2^{2k}(k!)^2}(x-\xi)^{2k}\nonumber\\&\hspace{10em}+\frac{2}{\pi}\sum_{k=1}^\infty(-1)^{k-1}\frac{K^{2k}}{2^{2k}(k!)^2}\Bigg(\sum_{j=1}^k\frac{1}{j}\Bigg)(x-\xi)^{2k}\nonumber\\&= \frac{2}{\pi}\ln(|x-\xi|)J_0(K|x-\xi|) + \sum_{k = 0}^{\infty}(-1)^k\frac{K^{2k}}{2^{2k-1}(k!)^2\pi}\Big(\ln(0.5K)+\gamma-h(k)\Big)
\end{align}
where
\begin{gather}
    h(k) = \begin{cases}
    0,&k = 0\\
    \displaystyle\sum_{j=1}^k\frac{1}{j},&k > 0.
    \end{cases}
\end{gather}
Setting
\begin{gather}\label{eq:fkK}
    f(k;K) = (-1)^k\frac{K^{2k}}{2^{2k-1}(k!)^2\pi}\Big(\ln(0.5K)+\gamma-h(k)\Big),
\end{gather}
we rewrite (\ref{eq:Y0}) in its final form
\begin{align}\label{eq:Y0-final}
    Y_0(K|x-\xi|) = \frac{2}{\pi}\ln(|x-\xi|)J_0(K|x-\xi|) +\sum_{k=0}^{\infty}f(K,k).
\end{align}
Substituting (\ref{eq:Y0-final}) into (\ref{eq:vossers-integral}),
\begin{align}\label{eq:vossers-integral-2}
    \int_{-L/2}^{L/2}\int_{-L/2}^{L/2}S''(x)S''(\xi)Y_0(K|x-\xi|)dxd\xi &= \frac{2}{\pi}\int_{-L/2}^{L/2}\int_{-L/2}^{L/2}S''(x)S''(\xi)\ln(|x-\xi|)J_0(K|x-\xi|)dxd\xi\nonumber\\
    &+ \sum_{k=0}^\infty f(k;K)I_k
\end{align}
with
\begin{align}\label{eq:Ik-0}
    I_k = \int_{-L/2}^{L/2}\int_{-L/2}^{L/2}S''(x)S''(\xi)(x-\xi)^{2k}dxd\xi.
\end{align}
We will now show $I_k$ can be expressed as a quadratic polynomial in moments of the type $M(n,0,0)$, in view of which we introduce the following notation,
\begin{align}\label{eq:moment_shortened_notation}
    M_n\equiv M(n,0,0).
\end{align}
To proceed, expand $(x-\xi)^{2k}$,
\begin{align}
    I_ k &= \int_{-L/2}^{L/2}\int_{-L/2}^{L/2}S''(x)S''(\xi)\sum_{i=0}^{2k}\binom{2k}{i}(-1)^ix^{2k-i}\xi^idxd\xi\nonumber\\
    &=\sum_{i=0}^{2k}\binom{2k}{i}(-1)^i\Bigg(\int_{-L/2}^{L/2}S''(x)x^{2k-i}dx\Bigg)\Bigg(\int_{-L/2}^{L/2}S''(\xi)\xi^id\xi\Bigg)\label{eq:Ik-1}.
\end{align}
It is evident that we must investigate the integral
\begin{gather}
    II_N = \int_{-L/2}^{L/2}S''(x)x^{N}dx.
\end{gather}
Integrating twice by parts,
\begin{align}
    II_N&= \Big[S'(x)x^{N}\Big]_{-L/2}^{L/2} -N\int_{-L/2}^{L/2}S'(x)x^{N-1}dx\nonumber\\
    &= \Big[S'(x)x^{N}\Big]_{-L/2}^{L/2} -N\Big[S(x)x^{N-1}\Big]_{-L/2}^{L/2} + N(N-1)\int_{-L/2}^{L/2}S(x)x^{N-2}dx\nonumber\\
    &=\Big[S'(x)x^{N}\Big]_{-L/2}^{L/2} -N\Big[S(x)x^{N-1}\Big]_{-L/2}^{L/2} + N(N-1)M_{N-2},\ N \ge 2\label{eq:IIN-1}.
\end{align}
Setting $M_{-1} \equiv M_{-2} = 0$ we extend (\ref{eq:IIN-1}) to all $N\ge 0$ without any effect on the ensuing calculations, which allows for less cumbersome notation in what follows. Now, looking at (\ref{eq:Ik-1}) it is evident that we must expand $II_{2k-i}II_i$,
\begin{align}
    II_{2k-i}II_{i} = &\Bigg(\Big[S'(x)x^{2k-i}\Big]_{-L/2}^{L/2} -(2k-i)\Big[S(x)x^{2k-i-1}\Big]_{-L/2}^{L/2} + (2k-i)(2k-i-1)M_{2k-i-2}\Bigg)\nonumber\\
    \cdot&\Bigg(\Big[S'(x)x^{i}\Big]_{-L/2}^{L/2} -i\Big[S(x)x^{i-1}\Big]_{-L/2}^{L/2} +i(i-1)M_{i-2}\Bigg)\nonumber\\
    = & \Bigg(\Big[S'(x)x^{2k-i}\Big]_{-L/2}^{L/2} -(2k-i)\Big[S(x)x^{2k-i-1}\Big]_{-L/2}^{L/2}\Bigg)\Bigg(\Big[S'(x)x^{i}\Big]_{-L/2}^{L/2} -i\Big[S(x)x^{i-1}\Big]_{-L/2}^{L/2}\Bigg)\nonumber\\[0.5em]
    + &(2k-i)(2k-i-1)\Bigg(\Big[S'(x)x^{i}\Big]_{-L/2}^{L/2} -i\Big[S(x)x^{i-1}\Big]_{-L/2}^{L/2}\Bigg)M_{2k-i-2}\nonumber\\[0.5em]
    +&i(i-1)\Bigg(\Big[S'(x)x^{2k-i}\Big]_{-L/2}^{L/2} -(2k-i)\Big[S(x)x^{2k-i-1}\Big]_{-L/2}^{L/2}\Bigg)M_{i-2}\nonumber\\[0.5em]
    +&i(i-1)(2k-i)(2k-i-1)M_{i-2}M_{2k-i-2}\label{eq:Ik-2}.
\end{align}
Before substituting (\ref{eq:Ik-2}) into (\ref{eq:Ik-1}), notice that a significant number of terms is repeated across the summation. For example, the last term in (\ref{eq:Ik-2}) at $i$ will be identical to that at $2k-i$. Same terms can easily be grouped together by noticing that for any $g(N)$
\begin{gather}\label{eq:sumoperator}
    \sum_{i=0}^{2k}\binom{2k}{i}(-1)^ig(i)g(2k-i) = \Bigg(\sum_{i=k}^k+2\sum_{i=0}^{k-1}\Bigg)\binom{2k}{i}(-1)^ig(i)g(2k-i).
\end{gather}
Then looking at (\ref{eq:Ik-2}) and (\ref{eq:sumoperator}), set
\begin{subequations}\label{eq:slenderbodycoefficients}
\begin{align}
    c^k &= \Bigg(\sum_{i=k}^k+2\sum_{i=0}^{k-1}\Bigg)\binom{2k}{i}(-1)^i\Bigg(\Big[S'(x)x^{2k-i}\Big]_{-L/2}^{L/2} -(2k-i)\Big[S(x)x^{2k-i-1}\Big]_{-L/2}^{L/2}\Bigg)\nonumber\\&\hspace{12em}\cdot\Bigg(\Big[S'(x)x^{i}\Big]_{-L/2}^{L/2} -i\Big[S(x)x^{i-1}\Big]_{-L/2}^{L/2}\Bigg),\\
    c^k_i &= \binom{2k}{i+2}(-1)^i(i+1)(i+2)\Bigg(\Big[S'(x)x^{2k-i-2}\Big]_{-L/2}^{L/2} -(2k-i-2)\Big[S(x)x^{2k-i-3}\Big]_{-L/2}^{L/2}\Bigg),\\
    c^k_{i,j} &= \binom{i+j+4}{i}(-1)^i(i+1)(i+2)(j+1)(j+2).
\end{align}
\end{subequations}
Transforming the operator of $I_k$ in (\ref{eq:Ik-1}) by (\ref{eq:sumoperator}) and substituting equations (\ref{eq:Ik-2}), (\ref{eq:slenderbodycoefficients}),
\begin{align}
    I_k = c^k + &\Bigg(\sum_{i=k}^k+2\sum_{i=0}^{k-1}\Bigg)\Bigg(c^k_{2k-i-2}M_{2k-i-2}\nonumber\\&\hspace{6em}+c^k_{i-2}M_{i-2}\nonumber\\&\hspace{6em}+c^k_{i-2,2k-i-2}M_{i-2}M_{2k-i-2}\Bigg)\nonumber\\
    = c^k + &\Bigg(\sum_{i=k}^k+2\sum_{i=0}^{k-1}\Bigg)\Bigg(c^k_{2k-i-2}M_{2k-i-2}\Bigg)\nonumber\\
    +&\Bigg(\sum_{i=k}^k+2\sum_{i=0}^{k-1}\Bigg)\Bigg(c^k_{i-2}M_{i-2}\Bigg)\nonumber\\
    +&\Bigg(\sum_{i=k}^k+2\sum_{i=0}^{k-1}\Bigg)\Bigg(c^k_{i-2,2k-i-2}M_{i-2}M_{2k-i-2}\Bigg)\nonumber\\
     = c^k + &c^k_{k-2}M_{k-2} + 2\sum_{i=0}^{k-1}\Bigg(c^k_{2k-i-2}M_{2k-i-2}\Bigg)\nonumber\\
     +&c^k_{k-2}M_{k-2} + 2\sum_{i=0}^{k-1}\Bigg(c^k_{i-2}M_{i-2}\Bigg)\nonumber\\
     +&c^{k}_{k-2,k-2}M_{k-2}^2 + 2\sum_{i=0}^{k-1}\Bigg(c^k_{i-2,2k-i-2}M_{i-2}M_{2k-i-2}\Bigg)\label{eq:Ik-3}.
\end{align}
To proceed, we reverse the order in first sum that appears in (\ref{eq:Ik-3}) and change the indices in the second and third sum by $i\to i+2$ which transforms the summation limits to $\{-2,k-3\}$.
\begin{align}
    I_k = c^k + &c^k_{k-2}M_{k-2} + 2\sum_{i=k-1}^{2k-2}\Bigg(c^k_{i}M_i\Bigg)\nonumber\\
     +&c^k_{k-2}M_{k-2} + 2\sum_{i=-2}^{k-3}\Bigg(c^k_iM_i\Bigg)\nonumber\\
     +&c^{k}_{k-2,k-2}M_{k-2}^2 + 2\sum_{i=-2}^{k-3}\Bigg(c^k_{i,2k-i-4}M_iM_{2k-i-4}\Bigg).
\end{align}
Notice that in the last two sums, the terms for $i=\{-2,-1\}$ originate from $II_0,\ II_1$ respectively and can be set to zero. Then,
\begin{align}
     I_k =\ &c^k + c^k_{k-2}M_{k-2} + 2\sum_{i=k-1}^{2k-2}\Bigg(c^k_{i}M_i\Bigg)\nonumber\\
     &\quad +c^k_{k-2}M_{k-2} + 2\sum_{i=0}^{k-3}\Bigg(c^k_iM_i\Bigg)\nonumber\\
     &\quad +c^{k}_{k-2,k-2}M_{k-2}^2 + 2\sum_{i=0}^{k-3}\Bigg(c^k_{i,2k-i-4}M_iM_{2k-i-4}\Bigg)\nonumber\\
     =\ &c^k\nonumber\\
     +&2\sum_{i=0}^{2k-2}\Bigg(c^k_iM_i\Bigg)\nonumber\\
     +&c^k_{k-2,k-2}M_{k-2}^2+2\sum_{i=0}^{k-3}\Bigg(c^k_{i,2k-i-4}M_iM_{2k-i-4}\Bigg)\label{eq:Ik-4}.
\end{align}
Finally, looking at (\ref{eq:vossers-integral-2}) and (\ref{eq:Ik-4}) we define the slender body operator ${\cal G}$ of order $n$ as,
\begin{gather}\label{eq:Gn}
    {\cal G}(n) = \sum_{k=0}^nf(k;K)I_k,
\end{gather}
which, under the assumption that the second derivative $S''(x)$  of the sectional area curve is bounded, will be a finite part of Vossers' integral. A proof for this statement is included in Appendix \ref{appendix:slenderbodyconvergentseries}

\subsection{Sobol's Sensitivity Analysis}\label{sec:SobolsSA}

\par\noindent As stated earlier, in this section Sobol's methodology for SA is presented, it's sampling-based counterpart which can be used to approximate the relevant sensitivity indices is provided and the approach through which parameters are categorised as sensitive or insensitive is introduced. Let $X_i$, $i=\{1,...,n\}$ be the random input variables to the function $f:R^n\to R^k$, which produces the random vector ${\bf Y} = f(X_1,...,X_{n})\in R^k$. Let ${\bf u}$ be an non-empty $r$-subset of $\{1,...,{n}\}$ and let ${\bf X_u} = (X_i,i\in {\bf u})\in R^r$, ${\bf X_{\tilde u}} = (X_i,i\in \{1,...,n\} / {\bf u}) \in R^{n-r}$. Then, recall the Hoeffding decomposition \cite{hoeffding1992class} of $f$:
\begin{gather}\label{eq:SobolHoeffingDecomp}
    f(X_0,...,X_{n-1}) = c + f_{\bf u}({\bf X_u}) + f_{\bf \tilde u}({\bf X_{\tilde u}}) + f_{\bf u,\tilde u}({\bf X_u}, {\bf X_{\tilde u}}),
\end{gather}
where $c\in R^k$, $f_{\bf u}:R^r\to R^k$, $f_{\bf \tilde u}:R^{n-r}\to R^k$ and $f_{\bf u,\tilde u}:R^n \to R^k$ are given by:
\begin{gather}
    c = E({\bf Y}),\quad f_{\bf u} = E({\bf Y}|X_{\bf u}) - c, \quad f_{\bf \tilde u} = E({\bf Y}|X_{\bf \tilde u}) - c,\quad f_{\bf u, \tilde u} = {\bf Y} - f_{\bf u}-f_{\bf u,\tilde u} - c.
\end{gather}
Taking the covariance of (\ref{eq:SobolHoeffingDecomp}), and due to $L^2$-orthogonality:
\begin{gather}
    Cov(f) = Cov(f_{\bf u}({\bf X_u})) + Cov(f_{\bf \tilde u}({\bf X_{\tilde u}})) + Cov(f_{\bf u,\tilde u}({\bf X_u}, {\bf X_{\tilde u}})) 
\end{gather}
Then to project these co-variances to a scalar, two approaches are to take the derivative or take the trace. The authors in \cite{GAMBOAMultiOutputs}, take the trace of these covariances and define the following three indices:
\begin{subequations}
\begin{equation}\label{eq:GenSobolIndexOnlyU}
    S^{\bf u}(f) = \frac{Trace(Cov(f_{\bf u}({\bf X_u})))}{Trace(Cov(f))},\quad \text{ sensitivity only to the inputs in ${\bf u}$}
\end{equation}
\begin{equation}\label{eq:GenSobolIndexOnlyNotU}
    S^{\bf \tilde u}(f) = \frac{Trace(Cov(f_{\bf \tilde u}({\bf X_{\tilde u}})))}{Trace(Cov(f))},\quad \text{ sensitivity only to the inputs not in ${\bf u}$}
\end{equation}
\begin{align}\label{eq:GenSobolIndexInteraction}
    S^{\bf u,\tilde u}(f) &= \frac{Trace(Cov(f_{\bf u,\tilde u}({\bf X_u}, {\bf X_{\tilde u}})))}{Trace(Cov(f))},\quad \text{  interaction between inputs of ${\bf u}$}\nonumber\\ &\hspace{15em}\text{and inputs of
$\{0, . . . , n-1\} / {\bf u}$}.
\end{align}
\end{subequations}
Now, for ${\bf u}=\{i\}$, (\ref{eq:GenSobolIndexOnlyU}) will be called the sensitivity index of the parameter $c_i$ with respect to the quantity of interest $f$, and denoted as,
\begin{gather}\label{eq:SIgeneraldef}
    SI_i^f \equiv S^{\bf u}(f).
\end{gather}
In accordance to \cite{GAMBOAMultiOutputs}, (\ref{eq:SIgeneraldef}) can be approximated via a sampling-based algorithm. Specifically, first draw $N$ independent samples ${\bf x}_i\in\mathbb{R}^n,\ i=1,...,N$ from the entire design space which are sampled according to the respective probability density function $\rho$ of the input parameters. Proceed by evaluating the respective responses ${\bf f}_i = f({\bf x}_i)$. Then, for each ${\bf x}_i$, construct another sample ${\bf x}_i^*$ by fixing $c_i$ and only varying the rest of the parameters to evaluate ${\bf f}_i^* = f({\bf x}_i^*)$. For ${\bf f}_i = (f_{i,1},...,f_{i,k})\in\mathbb{R}^k$, the sensitivity index $SI_i^f$ of $c_i$ is approximated by,
\begin{gather}\label{eq:SIapprox}
    SI_i^f = \frac{\displaystyle\sum_{j=1}^k\Bigg(\sum_{i=1}^Nf_{i,j}f_{i,j}^*-\frac{1}{N}\Bigg(\sum_{i=1}^N\frac{f_{i,j}+f_{i,j}^*}{2}\Bigg)^2\Bigg)}{\displaystyle\sum_{j=1}^k\Bigg(\sum_{i=1}^N\frac{(f_{i,j})^2+(f_{i,j}^*)^2}{2}-\frac{1}{N}\Bigg(\sum_{i=1}^N\frac{f_{i,j}+f_{i,j}^*}{2}\Bigg)^2\Bigg)}.
\end{gather}

\par\noindent Since in our case $c_1,\ c_2,\ c_3$ are design variables, $\rho$ is chosen as the uniform PDF \cite{Iooss2015} so as not to introduce any bias towards any region of the design space. As stated earlier, each parameter $c_i$ is characterised as sensitive or insensitive with respect to $f$ by  specifying a threshold value $\epsilon$ \cite{SARRAZIN2016GSA},
\begin{gather}\label{eq:threshold}
    \begin{cases}
        \text{SI}_i^f > \epsilon \implies c_i \text{ sensitive parameter with respect to $f$},\\
        \text{SI}_i^f \le \epsilon \implies c_i \text{ insensitive parameter with respect to $f$}.
    \end{cases}
\end{gather}
In alignment with \cite{GMGSA}, $\epsilon$ has been chosen as $0.05$ for this study.
\subsection{Physics-Geometry Correlation Measures}

\par\noindent Having evaluated the physics-based sensitivity indices $SI^{\cal P} = \{SI^{\cal P}_1,...,SI^{\cal P}_n\}$ and their geometry-based counterparts $SI^{\cal G} = \{SI^{\cal G}_1,...,SI^{\cal G}_n\}$ two metrics will be used to quantify their correlation,
\begin{subequations}
\begin{align}
    \text{NRMSE} &= \frac{\sqrt{\displaystyle\sum_{i=1}^n\frac{\Big(SI^{\cal P}_i-SI^{\cal G}_i\Big)^2}{n}}}{\max(SI^{\cal P})-\min(SI^{\cal P})}\label{eq:NRMSE},\\
    \text{similarity} &= \frac{\displaystyle\sum_{i=1}^n\overline{SI^{\cal P}_i}\cdot\overline{SI^{\cal G}_i}}{\displaystyle\sqrt{\sum_{i=1}^n\overline{SI^{\cal P}_i}}\sqrt{\sum_{i=1}^n\overline{SI^{\cal G}_i}}}\in[0,1]\label{eq:similarity},
\end{align}
where
\begin{gather}
    \overline{SI^{\cal P}_i} = \begin{cases}
        1,\ SI^{\cal P}_i > \epsilon\\ 
        0,\ \text{else}
    \end{cases},\quad \overline{SI^{\cal G}_i} = \begin{cases}
        1,\ SI^{\cal G}_i > \epsilon\\ 
        0,\ \text{else}.
    \end{cases}
\end{gather}
\end{subequations}
Equation (\ref{eq:NRMSE}) is the Root Mean Square Error between $SI^{\cal P}$ and $SI^{\cal G}$ normalised by the former, while (\ref{eq:similarity}) can be thought of as a discretised correlation between only the sensitive parameters of $SI^{\cal P}$ and $SI^{\cal G}$. More accurately, notice that if a parameter is insensitive in both $SI^{\cal P}$ and $SI^{\cal G}$ it then has no effect on (\ref{eq:similarity}). However if a parameter is sensitive for one of $SI^{\cal P}$, $SI^{\cal G}$ but not the other, (\ref{eq:similarity}) is penalised. This penalization can be interpreted as follows: if a parameter is sensitive for $SI^{\cal P}$ but not $SI^{\cal G}$ then the geometry-based DR has omitted important information regarding the relevant physics-based problem; if the opposite is true, then the geometry-based DR has not sufficiently reduced the dimensionality of the physics-based problem. These two metrics were used in \cite{GMGSA} in an analogous context to this study.

\section{Test Case: Wave Resistance and the Wigley-based Modeller}\label{sec:4}

\par\noindent Looking at Figure \ref{fig:MethodologyDiagram}, Section \ref{sec:3} was devoted to outlining two different choices of ${\cal G}$ as well as the approach to dimensionality reduction and correlation analysis between the physics- and geometry-based sensitivity indices. This section introduces the chosen parametric modeller ${\cal D}_1$ and physics-based property ${\cal P}$.

\subsection{Modified Wigley Hull Parametric Modeller}
\begin{figure}[h]
    \centering
    \includegraphics[width = 0.8\linewidth]{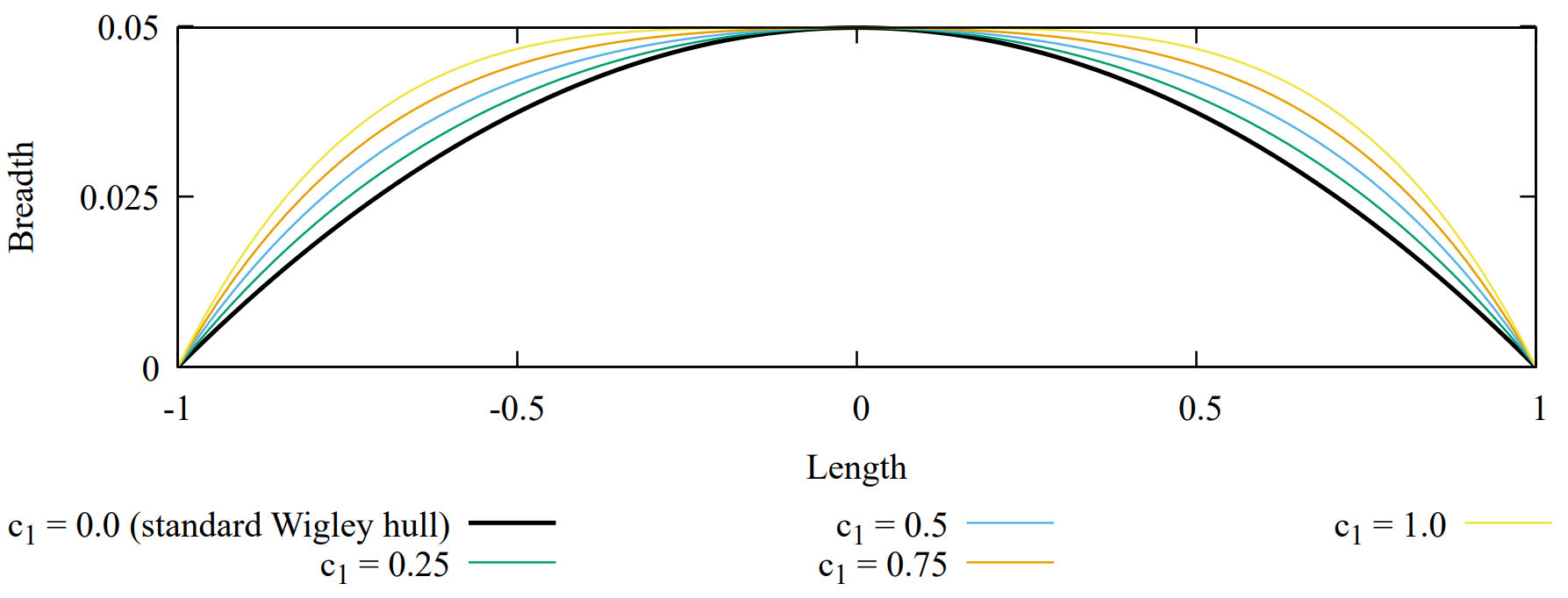}
    \caption{Effect of $c_1$ parameter on the modified Wigley hull. In this top-view figure, each curve is the deck-line corresponding to $c_1=\{0.0,\ 0.25,\ 0.5,\ 0.75,\ 1.0\}$ with all other parameters fixed at $\{L = 1,\ B = 0.0996,\ T = 0.13775,\ c_2=0,\ c_3=0\}$.}
    \label{fig:c1}
\end{figure}

\begin{figure}[h]
    \centering
    \includegraphics[width = 0.8\linewidth]{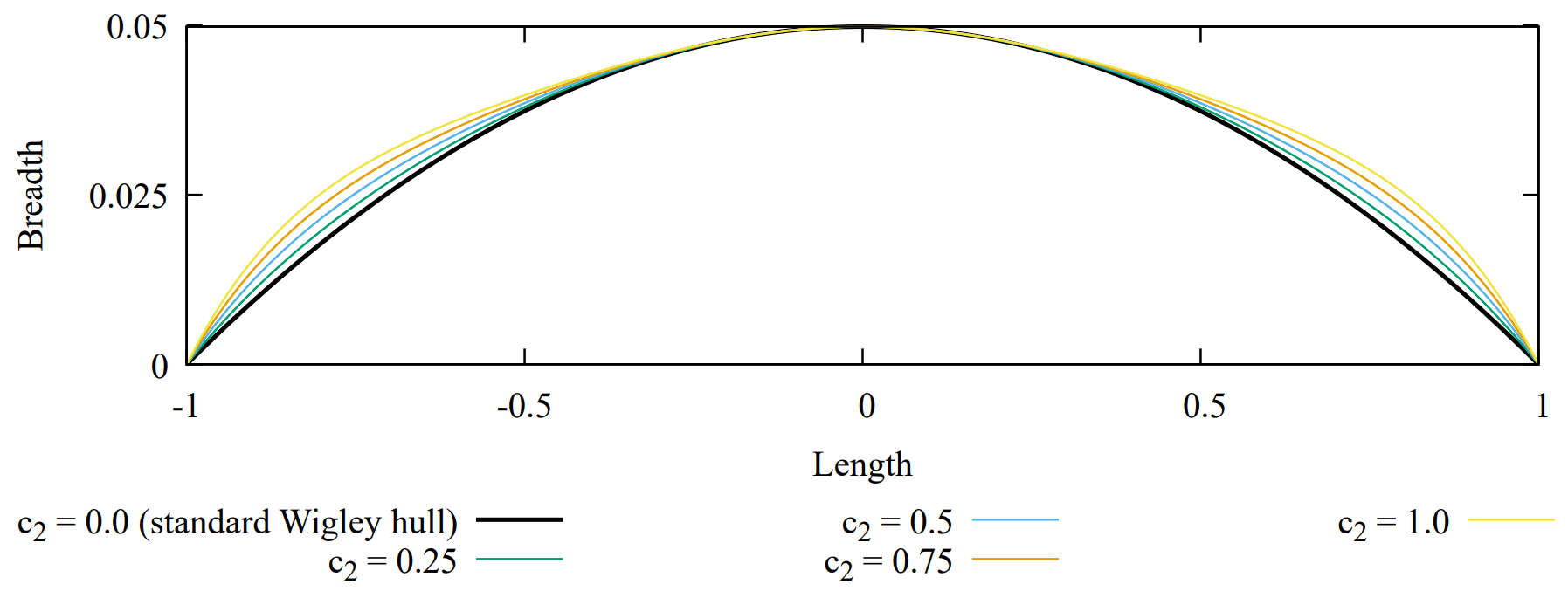}
    \caption{Effect of $c_2$ parameter on the modified Wigley hull. In this top-view figure, each curve is the deck-line corresponding to $c_2=\{0.0,\ 0.25,\ 0.5,\ 0.75,\ 1.0\}$ with all other parameters fixed at $\{L = 1,\ B = 0.0996,\ T = 0.13775,\ c_1=0,\ c_3=0\}$. }
    \label{fig:c2}
\end{figure}

\begin{figure}[h]
    \centering
    \includegraphics[width = 0.6\linewidth]{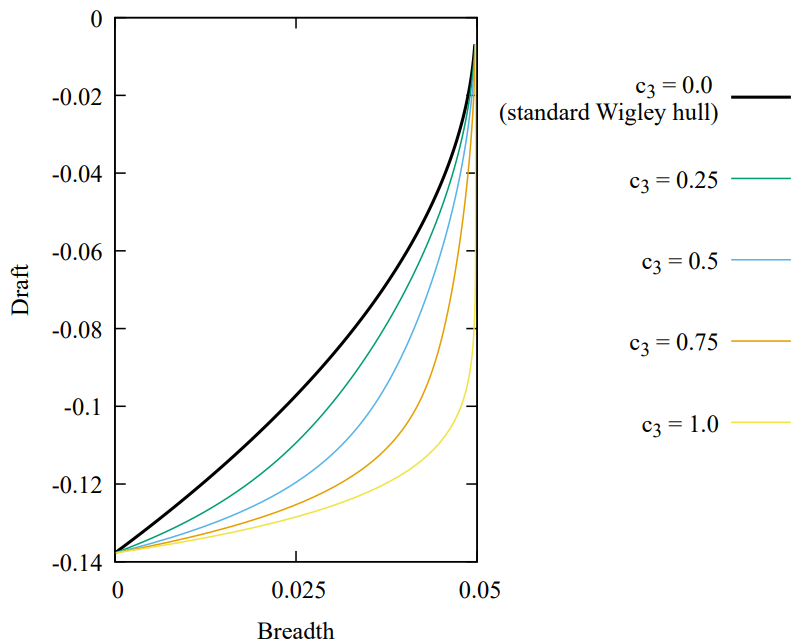}
    \caption{Effect of $c_3$ parameter on the modified Wigley hull. In this front-view figure, each curve is the depth-line corresponding to $c_3=\{0.0,\ 0.25,\ 0.5,\ 0.75,\ 1.0\}$ with all other parameters fixed at $\{L = 1,\ B = 0.0996,\ T = 0.13775,\ c_1=0,\ c_2=0\}$. }
    \label{fig:c3}
\end{figure}

The modeller ${\cal D}$ investigated in this study is a modified Wigley hull-form, a benchmark geometry which has been widely used is the literature \cite{journee1992experiments}. The half-length and half-breadth subsurface of $\cal D$, is parameterised in 6 parameters $\{L,B,T,c_1,c_2,c_3\}$ and given by,

\begin{subequations}\label{eq:wigley}
\begin{gather}
    {\cal D}(\xi,\zeta;L,B,T,c_1,c_2,c_3) = \begin{bmatrix}
    \displaystyle \frac{L\xi}{2}\\[1em]
    \displaystyle\frac{B\eta}{2}\\[1em]
    T\zeta
    \end{bmatrix}, \quad \xi\in[0,1],\ \zeta \in [0,1],\\
    \eta = (1-\zeta^2)(1-\xi^2)(1+c_1\xi^2+c_2\xi^4)+c_3\zeta^2(1-\zeta^8)(1-\xi^2)^4,
\end{gather}
\end{subequations}
where $L,\ B,\ T$ are the length breadth and draft respectively. While the effect of the first three parameters is clear, the same is not true for $c_1,\ c_2$ and $c_3$. Figures \ref{fig:c1}, \ref{fig:c2} and \ref{fig:c3} indicate this effect by varying each one between 0 and 1, while keeping the others fixed. Looking at the first two figures, note that $c_1$ and $c_2$ predominantly affect the deck area with negligible effect on the keel. By comparing Figures \ref{fig:c1} and \ref{fig:c2}, which are identically scaled, it is easy to see that $c_1$ has a stronger effect on the deck when compared to $c_2$. Varying $c_3$ mainly affects the keel of the hull, with no effect to the deck shape as is evident from the $\zeta^2$ factor in (\ref{eq:wigley}). In Figure \ref{fig:c3} one can see the transverse sections of the hull at midships ($\xi = 0$) for varying $c_3$. Finally, the original Wigley hull is recovered by setting $c_1=c_2=c_3=0$. 

It should be noted that if $S'(\pm L/2)\ne 0$, then (\ref{eq:vossers-integral-full}) diverges \cite{wehausen1973wave}, as is the case for the Wigley-hull family (\ref{eq:wigley}) where $S'(\pm L/2) = (8BT/3L)(1+c_1+c_2)\ne 0$ (see Appendix \ref{appendix:wigley_sac}). In this connection, we have adopted (\ref{eq:vossers-integral}) as a finite part of (\ref{eq:vossers-integral-full}).

\begin{figure}[h]
    \centering
    \includegraphics[width = 0.8\linewidth]{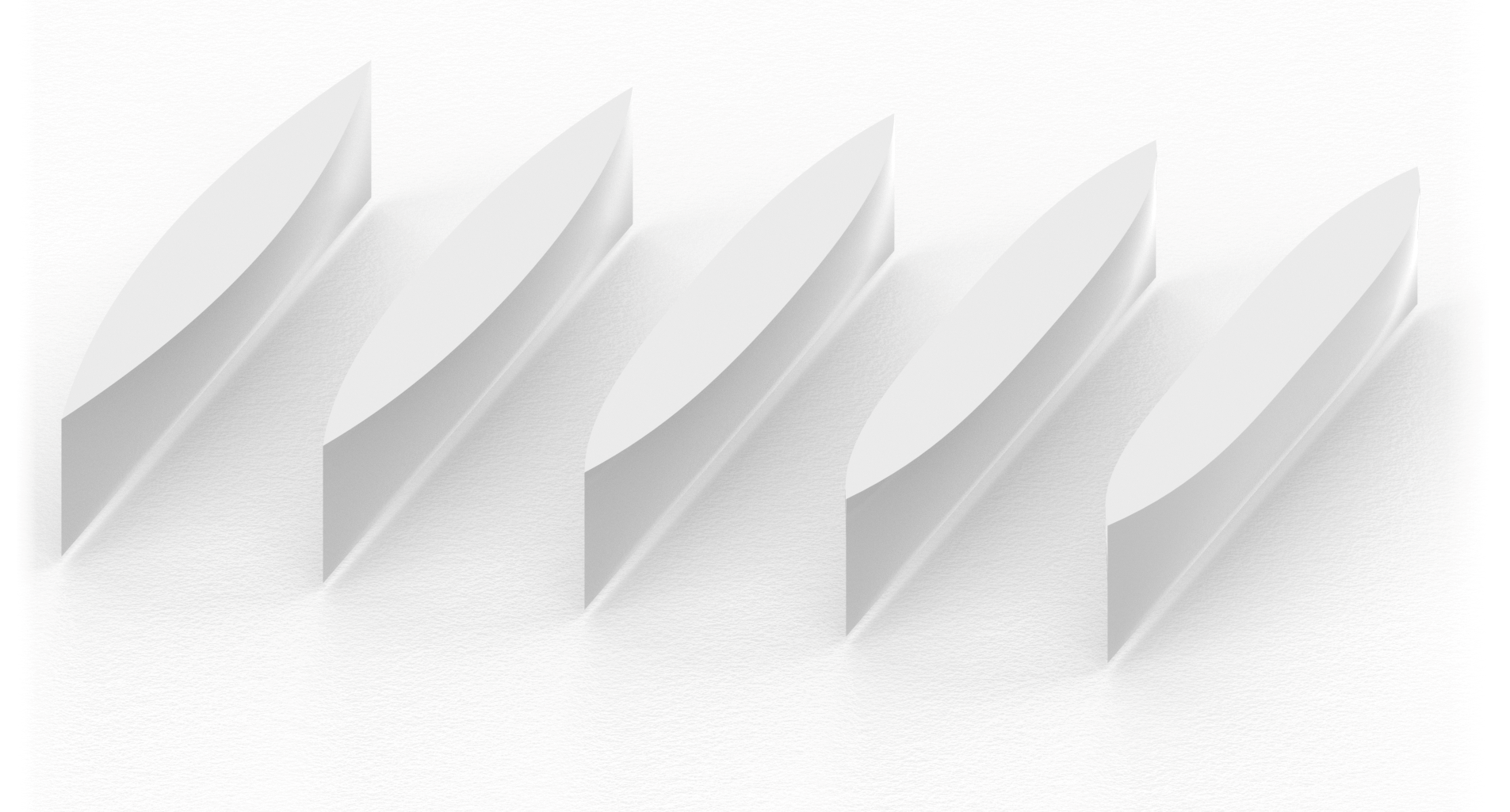}
    \caption{Renders of the modeller instances depicted in Figure \ref{fig:c1}. Specifically, from left to right, each hull corresponds to $c_1=\{0.0,\ 0.25,\ 0.5,\ 0.75,\ 1.0\}$ with all other parameters fixed at $\{L = 1,\ B = 0.0996,\ T = 0.13775,\ c_2=0,\ c_3=0\}$.}
    \label{fig:c1-render}
\end{figure}
\begin{figure}[h]
    \centering
    \includegraphics[width = 0.8\linewidth]{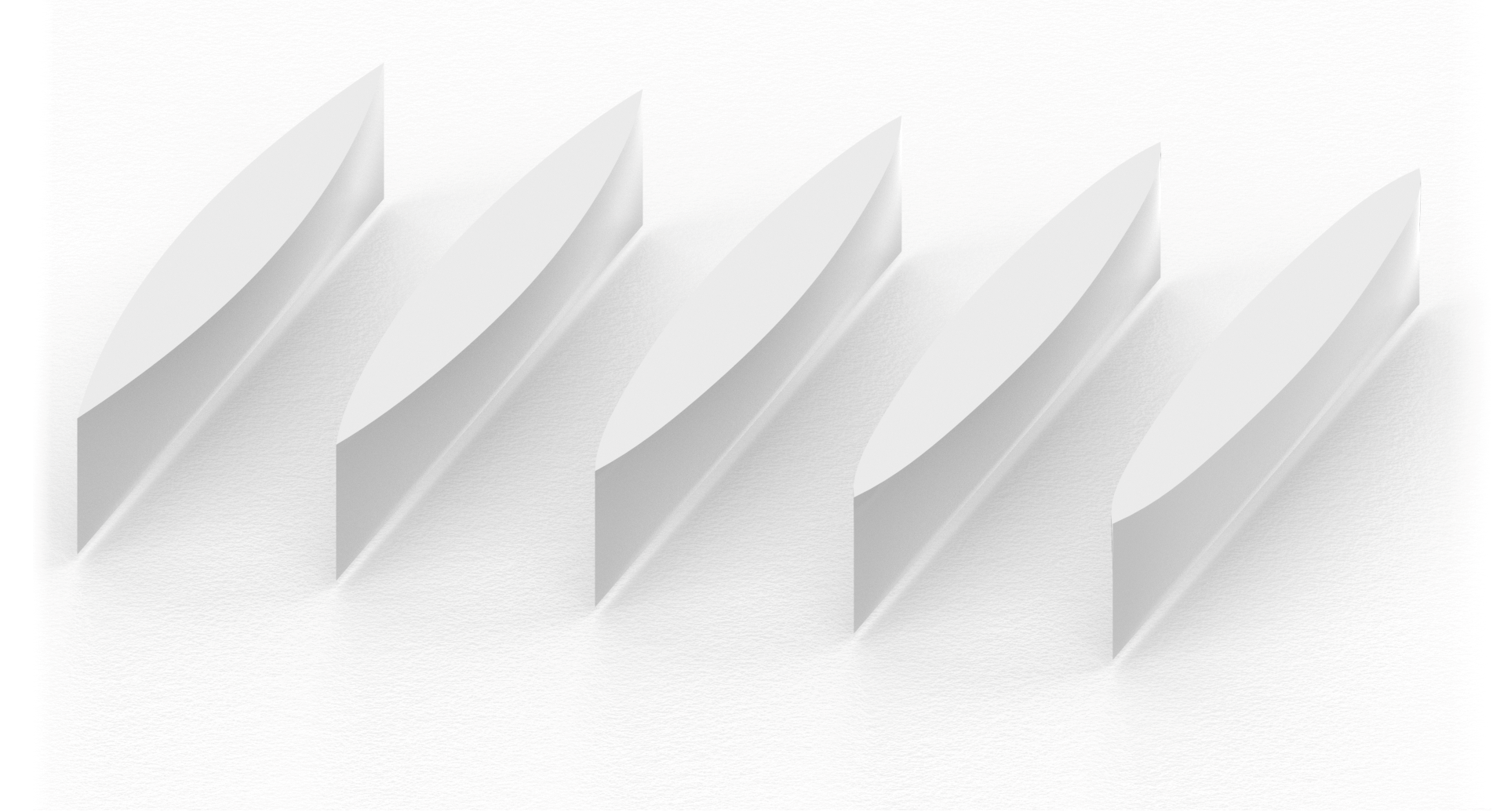}
    \caption{Renders of the modeller instances depicted in Figure \ref{fig:c2}. Specifically, from left to right, each hull corresponds to $c_2=\{0.0,\ 0.25,\ 0.5,\ 0.75,\ 1.0\}$ with all other parameters fixed at $\{L = 1,\ B = 0.0996,\ T = 0.13775,\ c_1=0,\ c_3=0\}$.}
    \label{fig:c2-render}
\end{figure}
\begin{figure}[h]
    \centering
    \includegraphics[width = 0.8\linewidth]{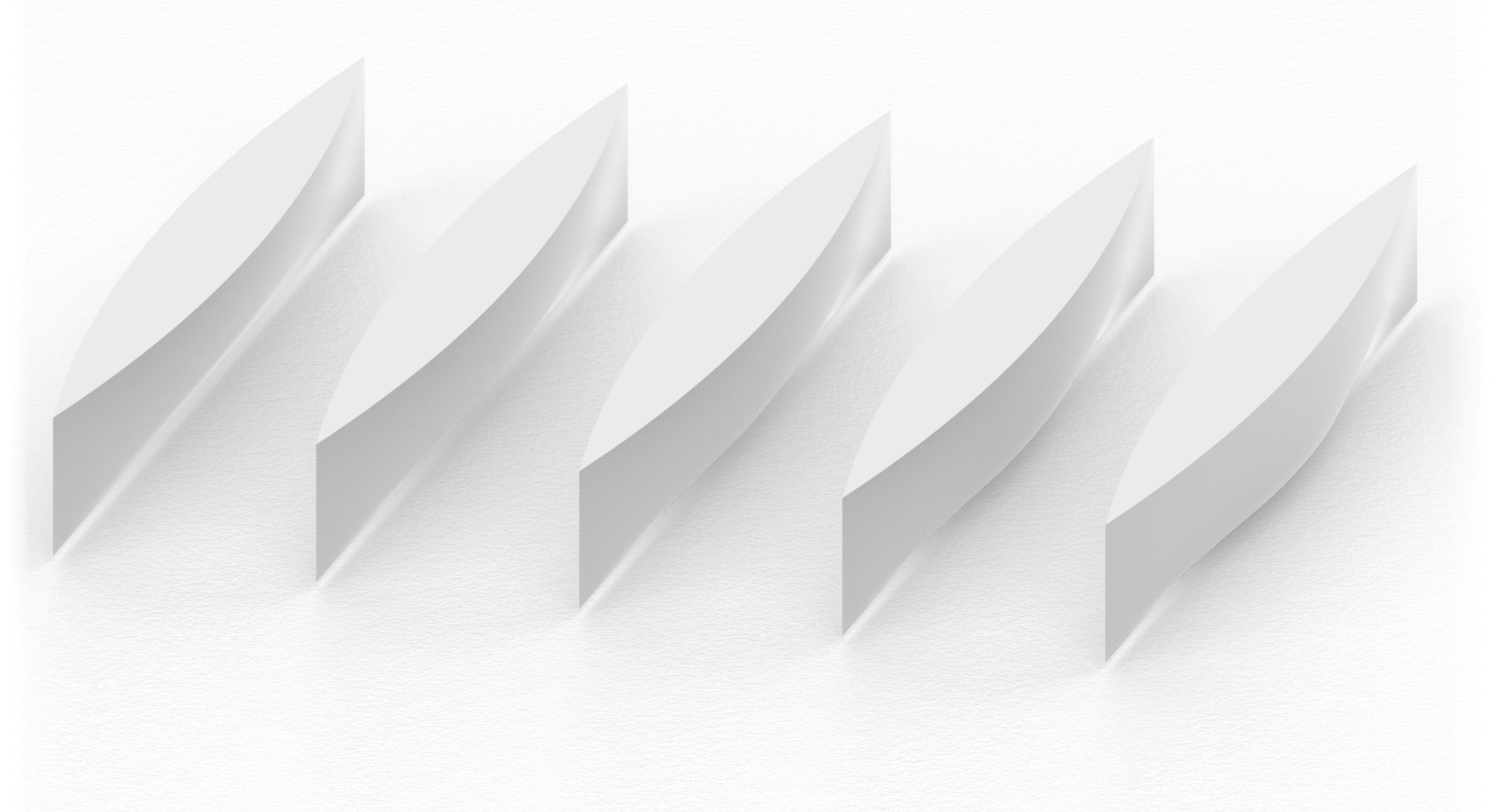}
    \caption{Renders of the modeller instances depicted in Figure \ref{fig:c3}. Specifically, from left to right, each hull corresponds to $c_3=\{0.0,\ 0.25,\ 0.5,\ 0.75,\ 1.0\}$ with all other parameters fixed at $\{L = 1,\ B = 0.0996,\ T = 0.13775,\ c_1=0,\ c_2=0\}$.}
    \label{fig:c3-render}
\end{figure}

\subsubsection{Design Space Selection}

\par\noindent The maximum number of parameters supported by (\ref{eq:wigley}) is 6, $\{L,\ B,\ T,\ c_1,\ c_2,\ c_3\}$. Similar to \cite{GMGSA}, a design space with uniform bounding box was chosen. We will now prove that this is indeed the case when $L,\ B$ and $T$ are fixed and $c_1\times c_2\times c_3 \in [0,1]^3$. \begin{subequations}
For fixed length, beam and draft, let ${\bf S}$ be the three-parameter family of surfaces defined by \begin{gather}
    {\bf S}(\xi,\zeta;c_1, c_2, c_3)\equiv {\cal D}(\xi,\zeta;L,B,T,c_1,c_2,c_3).
\end{gather} If $S_i$ is the $i^{th}$ coordinate function of ${\bf S}$, define the bounding box of ${\bf S}$ for a given choice of $c_1,\ c_2,\ c_3$ as
\begin{gather}
    {\cal B}(c_1,c_2,c_3) = \{{\bf p}=(p_1,p_2,p_3)^T\in\mathbb{R}^3:\ \min_{D_{\bf S}}S_i \le p_i \le \max_{D_{\bf S}}S_i,\quad i=1,2,3\},
\end{gather}
where $D_{\bf S}=[0,1]^2$ is the domain of ${\bf S}$. We will show that ${\cal B}(\cdot,\cdot,\cdot)$ is constant. Looking at (\ref{eq:wigley}), it is evident that for $i=1,3$ the minimum/maximum values of $S_i$ are taken at $\xi=\zeta=0$ / $\xi=\zeta=1$ respectively and are independent to $c_1,\ c_2,\ c_3$. For $i=2$ it is enough to show that $\displaystyle \min_{D_{\bf S}}\eta$ and $\displaystyle\max_{D_{\bf S}}\eta$ are independent to these three parameters as well. One can easily verify that in $D_{\bf S}$, $\eta \ge 0$ and $\eta(\xi,1;c_1,c_2,c_3)=0$. Therefore, $\displaystyle\min_{D_{\bf S}}\eta = 0$ which is independent to $c_1,\ c_2$ and $c_3$. For the maximum of $\eta$ notice that for $0\le c_i \le 1$,
\begin{gather*}
    0\le (1-\zeta^2)(1-\xi^2)(1+c_1\xi^2+c_2\xi^4) \le (1-\zeta^2)(1-\xi^2)(1+\xi^2+\xi^4)\\
    0\le c_3\zeta^2(1-\zeta^8)(1-\xi^2)^4\le \zeta^2(1-\zeta^8)(1-\xi^2)^4.
\end{gather*}
Adding the two inequalities yields
\begin{gather}\label{eq:wigleyboundingboxinvariancec1c2c3-dgdxi-2}
\eta(\xi,\zeta;c_1,c_2,c_3)\le \eta(\xi,\zeta;1,1,1),\quad \xi\times\zeta\times c_1\times c_2\times c_3 = [0,1]^5,
\end{gather}
where $\eta(\xi,\zeta;1,1,1)$ is clearly independent to $c_1,c_2$ and $c_3$. We complete this proof by showing that $\displaystyle\max_{D_{\bf S}}\eta(\xi,\zeta;c_1,c_2,c_3) = \max_{D_{\bf S}}\eta(\xi,\zeta;1,1,1)$ for any $c_1,\ c_2$ and $c_3$. Denote $g(\xi,\zeta)=\eta(\xi,\zeta;1,1,1)$ with $D_g=[0,1]^2$. Since $g$ is differentiable, if $g(\xi^*,\zeta^*)$ is a maximum, then either $(\xi^*,\zeta^*)$ is a singular point and/or it is on the boundary of $D_g$. We begin by investigating the first case,
\begin{align}
    \frac{\partial g}{\partial \xi} &= -2\xi(1-\zeta^2)(1+\xi^2+\xi^4)+(1-\zeta^2)(1-\xi^2)(2\xi+4\xi^3)+\zeta^2(1-\zeta^8)4(1-\xi^2)^3(-2\xi)\nonumber\\
    &= 2\xi(1-\zeta^2)\Big(-(1+\xi^2+\xi^4)+(1-\xi^2)(1+2\xi^2)-\zeta^2(1+\zeta^2)(1+\zeta^4)4(1-\xi^2)^3\Big)\nonumber\\
    &= 2\xi(1-\zeta^2)\Big(-1-\xi^2-\xi^4+1+2\xi^2-\xi^2-2\xi^4-\zeta^2(1+\zeta^2)(1+\zeta^4)4(1-\xi^2)^3\Big)\nonumber\\
    &=2\xi(1-\zeta^2)\Big(-3\xi^4-\zeta^2(1+\zeta^2)(1+\zeta^4)4(1-\xi^2)^3\Big)\label{eq:wigleyboundingboxinvariancec1c2c3-dgdxi-1}.
\end{align}
When not in the boundary of $D_g$ such that $\xi,\zeta\in(0,1)$ it is evident from (\ref{eq:wigleyboundingboxinvariancec1c2c3-dgdxi-1}) that ${\partial g}/{\partial \xi} < 0$ and therefore, there are no singular points in the interior of $D_g$. For the second case where $(\xi,\zeta)\in{\rm Bd}(D_g)$,
\begin{gather}
    g({\rm Bd}(D_g)) = \begin{Bmatrix}0&,&\quad (\xi,\zeta)=(1,\zeta)\ {\rm or}\ (\xi,1)\\1-\zeta^8&,&(\xi,\zeta)=(0,\zeta)\\
    1-\xi^6&,&(\xi,\zeta)=(\xi,0)\end{Bmatrix}.
\end{gather}
Then $\displaystyle\max_{{\rm Bd}(D_g)}g = 1$ and further, by the non-singularity of $g$ in ${\rm Int}(D_g)$, it follows that $\displaystyle\max_{D_g}g=\max_{{\rm Bd}(D_g)}g=1$. All that is left now is to show that $\displaystyle\max_{D_{\bf S}}\eta(\xi,\zeta;c_1,c_2,c_3) = \max_{D_{\bf S}}\eta(\xi,\zeta;1,1,1) {\buildrel {\rm def }\over=} \max_{D_g}g(\xi,\zeta) = 1$ for any $c_1,\ c_2$ and $c_3$. By (\ref{eq:wigleyboundingboxinvariancec1c2c3-dgdxi-2}) it is sufficient to show that there exist $\xi^*$ and $\zeta^*$ such that $\eta(\xi^*,\zeta^*;c_1,c_2,c_3) = 1$. Indeed, substitute $(\xi^*,\zeta^*)=(0,1)$ in equation (\ref{eq:wigley}) to verify this statement and complete the proof.
\end{subequations}

\subsubsection{Closed Forms of the Geometric Moments of the Modified Wigley Hull}

\par\noindent Due to the simple analytical expression (\ref{eq:wigley}) of the modified Wigley hull, it is possible to evaluate the moments (\ref{eq:MIpqr}) analytically with the relevant analysis carried out in Appendix \ref{appendix:closedformsappendix}. Specifically, $M(p,q,r)$ and $M_T(p,q,r)$ are given in equations (\ref{eq:Mpqr-final-FULL}) and (\ref{eq:MTpqr-final-FULL}) respectively, while $M_S(p,q,r)$ and $M_I(p,q,r)$ can easily be computed using (\ref{eq:MSpqr}) and (\ref{eq:MIpqr})

\subsection{Wave Resistance Coefficient}

\par\noindent In line with \cite{GMGSA}, ${\cal P}$ is identified with the wave resistance coefficient, $C_w$. This choice can be attributed to its explicit connection to the SAC, $S(x)$, of a hull which is a purely geometric property. This connection was initially established in the context of slender body theory via the so called Vossers' integral \cite{vossers1962some} as in equation (\ref{eq:vossers-integral-full}). Further development continued (\cite{maruo1962calculation}; \cite{tuck1963vossers}; \cite{tuck1964systematic}) which, however, resulted in a theory suffering from a number of deficiencies as pointed out in \cite{kotik1963various} and \cite{wehausen1973wave}.  Nevertheless, the SAC's importance in resistance-related ship-design problems which originated by the so-called Lackenby transformation \cite{lackenby1950systematic} has continued till recent times (\cite{han2012hydrodynamic}; \cite{tasrief2013improvementofshipgeometry})\\

\par\noindent The $C_w$ of (\ref{eq:wigley}) was approximated using a Boundary Element Method-based solver in the context of IsoGeometric Analysis (IGA-BEM), as introduced in \cite{BELIBASSAKIS201353}. To utilise this solver, the input geometry must be provided in B-spline form, which is not the case for (\ref{eq:wigley}). However, it is important to note that (\ref{eq:wigley}) represents a polynomial surface of degree 8 in $\xi$ and 10 in $\zeta$. This means that only a change from the monomial to the B-spline basis is necessary or rather to the simpler Bernstein basis, needing only a single $8\times10$ patch for the surface to be explicitly represented. The resulting Bezier surface will be represented by $(8+1)(10+1)=99$ control points ${\bf b}_{i,j}$, $i=0,...,8,\ j=0,...,10$ which are calculated numerically: for $B_{n,i}$ the Bernstein polynomial of degree n, find ${\bf b}_{i,j}$ such that for ${\cal D}:[0,1]^2\to\mathbb{R}^3$ the wigley hull presented in (\ref{eq:wigley}) it is true that
\begin{gather}\label{eq:BsplineExtraction-1}
    {\cal D}(\xi,\zeta) = \sum_{i=0}^8\sum_{j=0}^{10}B_{8,i}(\xi)B_{10,j}(\zeta){\bf b}_{i,j},\quad \xi\times\zeta \in [0,1]^2.
\end{gather}
To find the control points in (\ref{eq:BsplineExtraction-1}) it is enough to evaluate $\cal D$ at 99 distinct points in $[0,1]^2$, ${\bf p}_{i,j}$ so that (\ref{eq:BsplineExtraction-1}) is converted to three linear systems (one for each coordinate of ${\bf b}_{i,j}$) each of 99 independent equations in 99 unknowns. Finally, the choice of ${\bf p}_{i,j}$ was done uniformly in $\xi,\ \zeta$,
\begin{gather}\label{eq:BsplineExtraction-2}
    {\bf p}_{i,j} = (i/8, j/10),\quad i=0,...,8,\ j=0,...,10.
\end{gather}
Having found ${\bf b}_{i,j}$ it was then straightforward to supply the relevant knot sequences such that the resulting Bezier surface can be interpreted as a B-spline patch.

\section{Results and Discussion}\label{sec:5}

\par\noindent This section presents the detailed numerical results of the methodology outlined in Figure \ref{fig:MethodologyDiagram} both for SSV- and slender body (${\cal G}$)-based geometric operators, presented in Sections \ref{sec:SSV} and \ref{sec:SlenderBodyOperator}, respectively. The approximation of the relevant sensitivity indices is carried out using a sampling-based strategy. This implies the necessity of selecting a sampling strategy and analysing the computed approximation (\ref{eq:SIapprox}) for convergence.

The chosen sampling strategy is detailed in Section \ref{sec:DPSsection}, while its convergence in case of wave resistance, SSV- and {\cal G}-based geometric operators is presented in Sections \ref{sec:CwBasedResults}, \ref{sec:SSVbasedResults}, and \ref{sec:SBOresults}, respectively. Finally, in Section \ref{sec:summaryOfSAresults}, all the aforementioned results are summarised. Correlation measures (\ref{eq:NRMSE}) and (\ref{eq:similarity}) are applied to the computed indices to compare the performance of SSV versus that of the slender body operator.

\subsection{Dynamic Propagation Sampling}\label{sec:DPSsection}

The choice of sampling algorithm is of fundamental importance to the expedited convergence of (\ref{eq:SIapprox}), therefore, a Dynamic Propagation Sampling (DPS) \cite{KHANintrasensitivity} is utilised in this study. The principal idea behind DPS is an optimisation based approach, which aims at creating uniformly distributed and diverse number of sample over multiple run. DPS achieves this while minimising an objective function, which is a weighted combination of three sampling criteria; \textbf{space-filling}, \textbf{ non-collapsing} and \textbf{repulsion}.

Given a sample-set of $N$ samples, ${\bf S}_N = \{{\bf s}_i,\ i=1,...,N\}$, ${\bf s}_i\in\mathbb{R}^n$ the space-filling criterion ensures that the samples ${\bf s}_i$ are uniformly distributed over the entire design space. This criterion is defined as $f_1:(\mathbb{R}^n)^N\to\mathbb{R}$
\begin{align}
    f_1({\bf S}_N) = \sum_{i=1}^{N-1}\sum_{j=i+1}^N\frac{1}{\|{\bf s}_i-{\bf s}_j\|^2},
\end{align}
where $\|\cdot\|$ is the Euclidean norm in $\mathbb{R}^k$. However, minimisation of $f_1$ tends to place samples towards the boundary of the design space. To address this issue, the non-collapsing criterion is introduced through the discretisation of the design space into $N^k$ disjoint subsets and the construction of a functional that discourages samples from being in the same subset. Specifically, if the range of each parameter $c_i\in[a_i,b_i]$ is divided into $N$ sub-intervals $a_i = t_i^0 < t_i^1 < ... < t_i^N = b_i$ define the discrete position ${\bf d}_s\in\mathbb{R}^n$ of a sample ${\bf s}\in\mathbb{R}^n$ as,
\begin{subequations}
\begin{gather}
    {\bf d}_{\bf s} = [d_1,...,d_N],\quad \text{such that } \forall i,\quad {\bf s}_i \in [t_i^{d_i},t_i^{d_i+1}).
\end{gather}
Then, the related functional is defined as,
\begin{align}
    f_2({\bf S}_N) &= \omega \sum_{i=1}^{N-1}\sum_{j=i+1}^N{\cal K}({\bf s}_i,{\bf s}_j),\\
    {\cal K}({\bf s}_i,{\bf s}_j) &= \sum_{m=1}^n \delta\Big(\big({\bf d}_{{\bf s}_i}\big)_m,\big({\bf d}_{{\bf s}_j}\big)_m\Big),
\end{align}
\end{subequations}
where $\omega$ is a user-controlled weight on the non-collapsing criterion, $\big({\bf d}_{\bf s}\big)_m$ is the $m^{th}$ element of the discrete position of ${\bf s}$ and $\delta(i,j)$ is Kronecker's delta. 

Finally, DPS can be used to create a new set of samples, while ensuring that they are different from the previous sampled samples. Given an already existent sample set of $M$ samples $\hat{\bf S}_M=\{\hat{\bf s}_i,\ i=1,...,M\}$ DPS uses the repulsion criterion ($f_3$) to ensure that the new samples ${\bf S}_N$ are different from the previous sampled designs. The formulation of $f_3$ is similar to that of $f_1$,
\begin{align}
    f_3({\bf S}_N) = \sum_{i=1}^{N}\sum_{j=1}^M\frac{1}{\|{\bf s}_i-\hat{\bf s}_j\|^2},
\end{align}
where $\|\cdot\|$ is the Euclidean norm in $\mathbb{R}^k$. Finally, define $f({\bf S}_N) = f_1({\bf S}_N) + f_2({\bf S}_N) + f_3({\bf S}_N)$ as the functional to be minimised.\\

\subsection{$C_w$ based SA}\label{sec:CwBasedResults}

\par
\begin{figure}[H]
    \centering
    \includegraphics{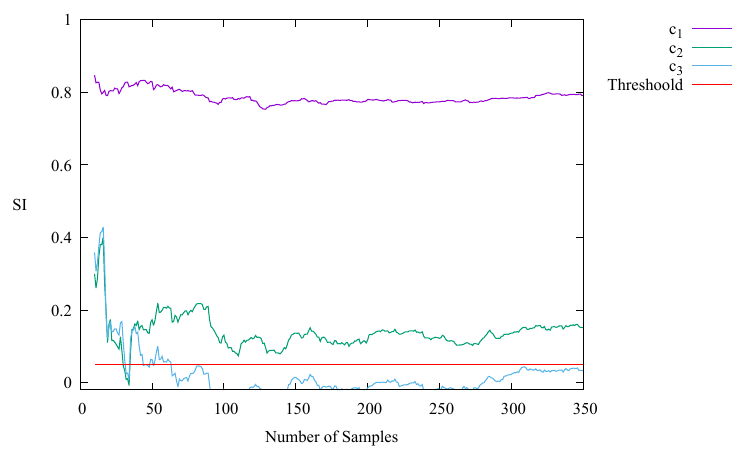}
    \caption{Convergence of $C_w$-based sensitivity indices. The horizontal axis is the number of DPS-generated samples while the vertical axis shows the sensitivity indices (SI)}
    \label{fig:CwConvergence}
\end{figure}

\noindent In order to identify the number of samples $N$ to use for the reliable implementation of (\ref{eq:SIapprox}) the relevant convergence has to be investigated. Figure \ref{fig:CwConvergence} show the sensitivity indices of the parameters with respect to $C_w$ over the varying number of samples. Evaluation of sensitivity indices commenced with 10 samples and varied to 350 samples. It can be seen that sensitivity indices are unstable until 100 samples, however, they tend to converge afterwarding, indicating that $c_1,\ c_2$ are sensitive parameters while $c_3$ is an insensitive parameter. This is to be expected due to the nature of the wave-making phenomenon: due to the wave effect exponential decaying with respect to depth, $c_3$ which as shown in Figure (\ref{fig:c3}) predominantly affects the keel of the hull, is expected to have less of an effect on $C_w$ when compared to $c_1,\ c_2$ which predominantly affect the deck. Moreover, notice that the effect of $c_1$ on $C_w$ is significantly greater than that of $c_2$ which can be again justified by comparing Figures \ref{fig:c1} and \ref{fig:c2}: the $c_1$ parameter has a significantly greater effect on the deck-geometry, compared to $c_2$.

\subsection{SSV based SA}\label{sec:SSVbasedResults}

\par
\begin{figure}[H]
    \centering
    \includegraphics{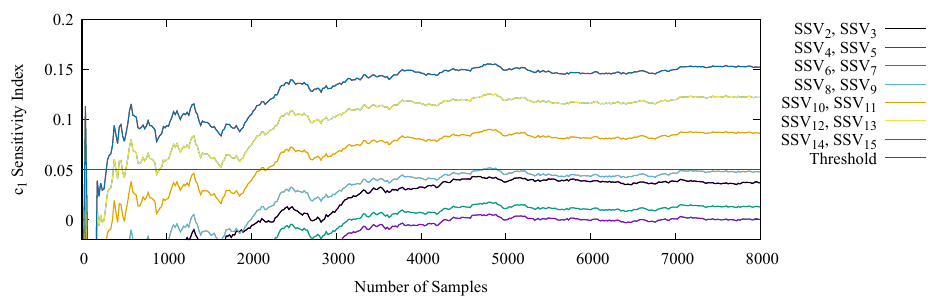}
    \caption{Convergence of the SSV-based sensitivity index of parameter $c_1$ with respect to the number of samples across varying order from 2 to 15.}
    \label{fig:SSVc1Convergence}
\end{figure}

\begin{figure}[H]
    \centering
    \includegraphics{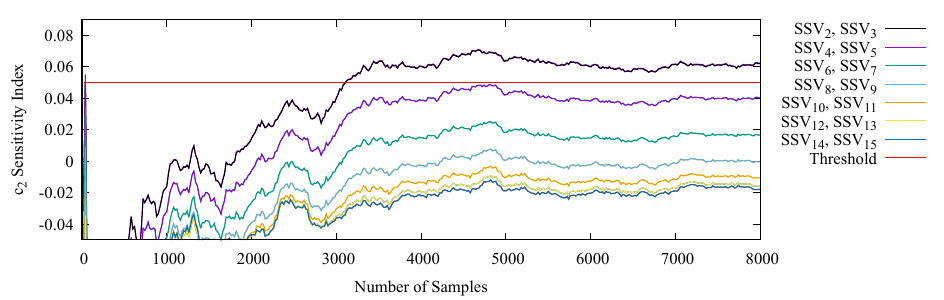}
    \caption{Convergence of the SSV-based sensitivity index of parameter $c_2$ with respect to the number of samples across varying order from 2 to 15.}
    \label{fig:SSVc2Convergence}
\end{figure}

\begin{figure}[H]
    \centering
    \includegraphics{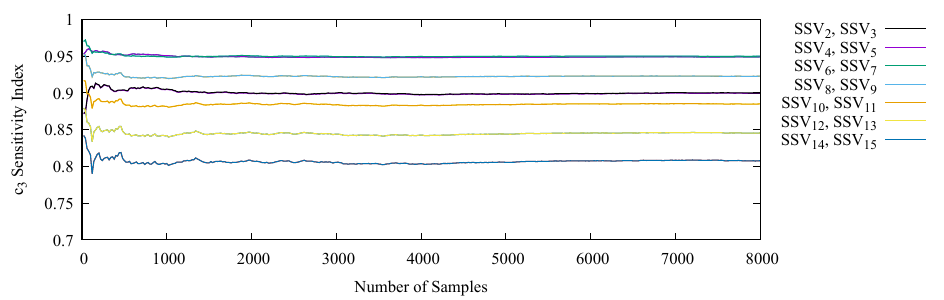}
    \caption{Convergence of the SSV-based sensitivity index of parameter $c_3$ with respect to the number of samples across varying order from 2 to 15.}
    \label{fig:SSVc3Convergence}
\end{figure}

\noindent In contrast to the $C_w$-based SA, performing sensitivity analysis with respect to SSV necessitates the choice of a right order $n$. In \cite{GMGSA}, the authors investigated up to $n=4$, as high-order geometric moments can be sensitive to noise while at the same time, numerical inaccuracies are ever-present when evaluating high-order terms. However, due to the simple expression of (\ref{eq:wigley}), the evaluation of high-order moments is computationally inexpensive and less pronounce to noise and inaccuracies. Therefore, in this study, the chosen range of $n$ is $\{2,\cdots,15\}$. Consequently, each order $n$ will yield its respective parameter sensitivity indices, necessitating a convergence analysis for all orders and parameters. This results in a total of 42 distinct convergence graphs, which can be found in Figures \ref{fig:SSVc1Convergence} through \ref{fig:SSVc3Convergence}. Each graph corresponds to one of the three parameters: $c_1$, $c_2$, and $c_3$.

Several noteworthy observations should be made about these results. Firstly, it is apparent that the number of samples required for convergence is significantly higher compared to $C_w$ in Figure \ref{fig:CwConvergence}. This disparity can be attributed to certain sensitivity indices being closer to the threshold value set at 0.05. Additionally, indices evaluated with respect to higher-order SSVs tend to exhibit greater stability with regard to the number of samples. In other words, they require fewer samples to achieve convergence compared to lower-order SSVs. This trend is particularly evident in parameters $c_1$ and $c_2$, which were sensitive to $C_w$. Moreover, note the similarity in sensitivity indices related to consecutive orders of the form $2k$ and $2k+1$. This behaviour can be attributed to the length and breadth-wise symmetry of the Wigley hull (\ref{eq:wigley}), resulting in the vanishing of odd-index moments in their first and/or second positions, as seen in (\ref{eq:Mpqr-final-FULL}). 

\par
\begin{figure}[H]
    \centering
    \includegraphics{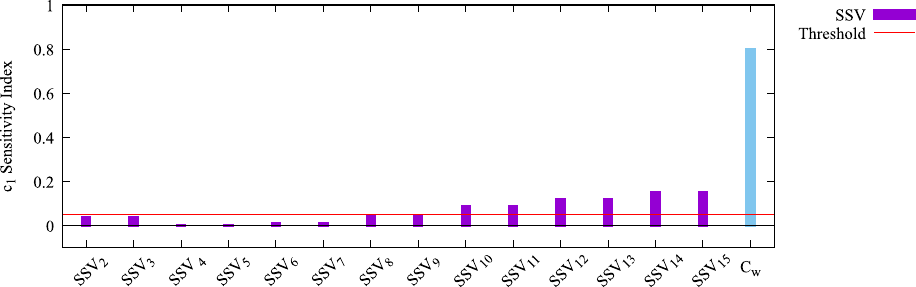}
    \caption{$C_w$- and SSV-based sensitivity indices of parameter $c_2$ evaluated with 8000 samples over SSV orders varying from 2 to 15.}
    \label{fig:SSVc1}
\end{figure}

\begin{figure}[h!]
    \centering
    \includegraphics{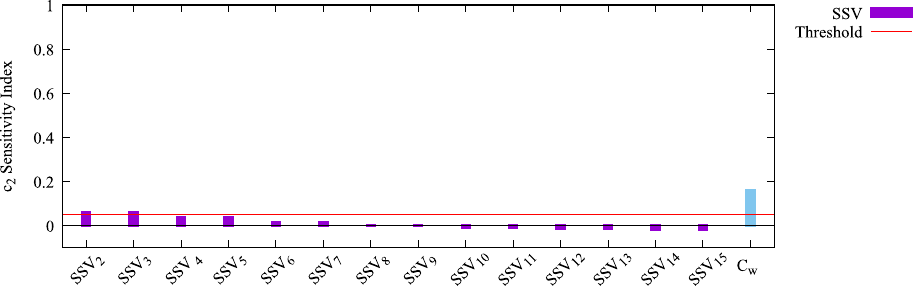}
    \caption{$C_w$- and SSV-based sensitivity indices of parameter $c_2$ evaluated with 8000 samples over SSV orders varying from 2 to 15.}
    \label{fig:SSVc2}
\end{figure}

\begin{figure}[h!]
    \centering
    \includegraphics{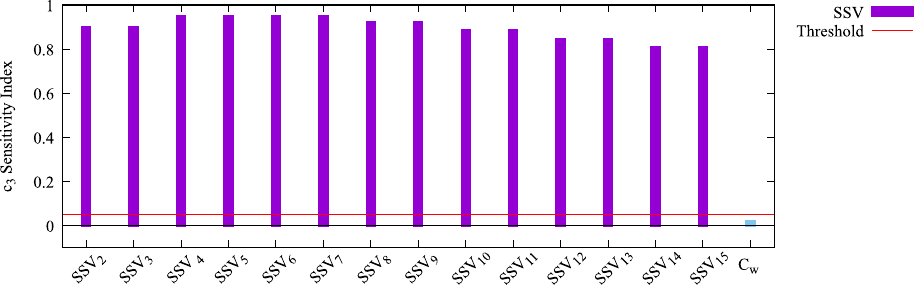}
    \caption{$C_w$- and SSV-based sensitivity indices of parameter $c_2$ evaluated with 8000 samples over SSV orders varying from 2 to 15.}
    \label{fig:SSVc3}
\end{figure}

\noindent Figures \ref{fig:SSVc1} through \ref{fig:SSVc3} show the comparative results of sensitivity indices with respect to SSV and $C_w$, where each figure is related to one of the parameters. In the case of parameter $c_1$, which is clearly sensitive to $C_W$, its sensitivity index increases with the order $n$ of SSV. Interestingly, it is insensitive until $\text{SSV}_7$ and then tends to increase steadily. However, in Figure \ref{fig:SSVc2}, $c_2$ indicates a behaviour opposite to $c_1$, i.e., its sensitivity reduces as the order of SSV increases. Based on the threshold of 0.05, it tends to be sensitive up to $\text{SSV}_5$. Note that in Figure \ref{fig:SSVc2}, some of the indices are negative, indicating convergence, but some are numerically zero. Since their confidence intervals include zero, we can safely treat them as zero. Finally, in the case of $c_3$, SSV$_n$ and $C_w$ are in complete disagreement, with this parameter being insensitive to $C_w$ while highly sensitive to SSV$_n$. This misalignment in the $c_3$ highlights one of the shortcomings of SSV—its inability to differentiate between keel-wise and deck-wise changes to the geometry, which are crucial for wave-making.

\subsection{Slender Body Operator based SA}\label{sec:SBOresults}

\begin{figure}[H]
    \centering
    \includegraphics{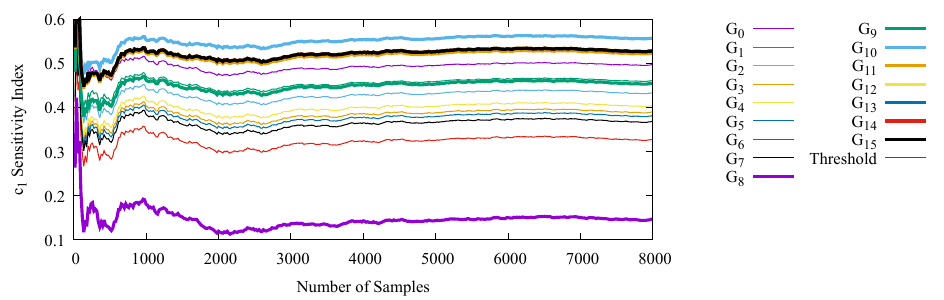}
    \caption{Convergence of the ${\cal G}$-based sensitivity index of parameter $c_1$ with respect to the number of samples across different orders of ${\cal G}$.}
    \label{fig:Gc1Convergence}
\end{figure}

\begin{figure}[H]
    \centering
    \includegraphics{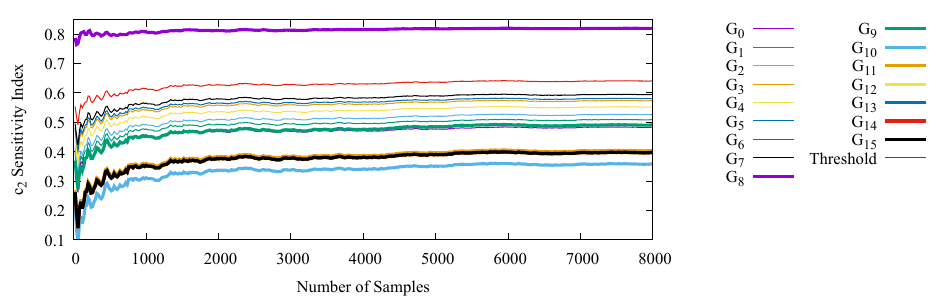}
    \caption{Convergence of the ${\cal G}$-based sensitivity index of parameter $c_2$ with respect to the number of samples across different orders of ${\cal G}$.}
    \label{fig:Gc2Convergence}
\end{figure}

\begin{figure}[H]
    \centering
    \includegraphics{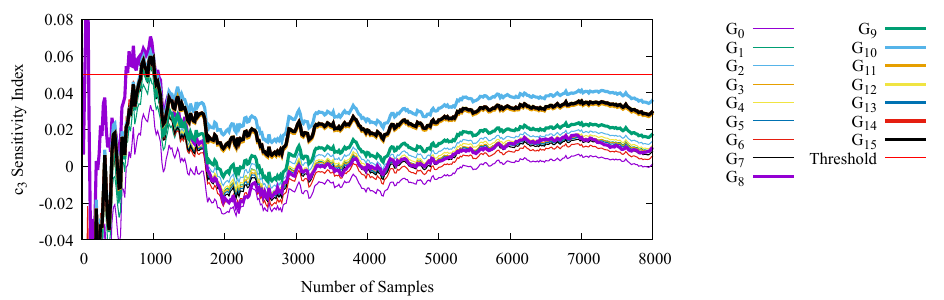}
    \caption{Convergence of the ${\cal G}$-based sensitivity index of parameter $c_3$ with respect to the number of samples across different orders of ${\cal G}$.}
    \label{fig:Gc3Convergence}
\end{figure}

Considering (\ref{eq:Gn-convergence-6}) and in contrast to the SSV-based sensitivity analysis, it would not yield meaningful results to the reader by deriving outcomes for a range of orders $n$ as depicted in Figures \ref{fig:SSVc1Convergence} through \ref{fig:SSVc3}. This would unnecessarily lengthen the paper. Instead, it is more reasonable to identify a singular sufficiently high $n$ where the corresponding sensitivity indices have reached convergence. To ensure the convergence of (\ref{eq:SIapprox}) concerning both the sample size and the order $n$. These two variables are incrementally adjusted in the following manner: for each order $n$, the necessary sample size to guarantee the convergence of ${\cal G}(n)$ is determined before examining whether the consecutive differences $|{\cal G}(k+1)-{\cal G}(k)|$ have sufficiently diminished. Figures \ref{fig:Gc1Convergence} through \ref{fig:Gc3Convergence} depict the ${\cal G}$-based sensitivity indices $c_1$, $c_2$, and $c_3$ respectively, as the sample size increases for ${\cal G}$-orders ranging from 0 to 15.

\begin{figure}[H]
    \centering
    \includegraphics{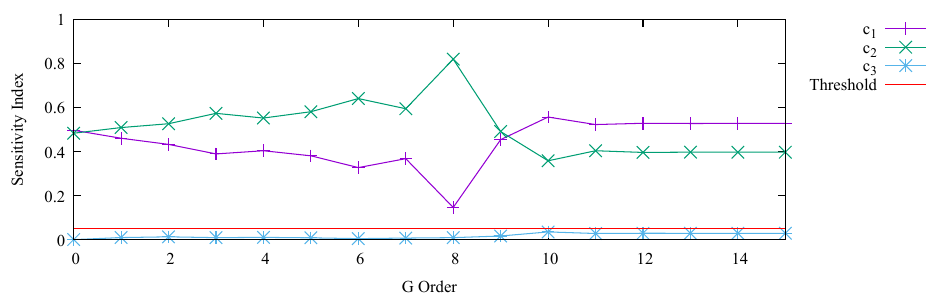}
    \caption{Convergence of ${\cal G}$-based sensitivity indices with respect to order from 0 to 15 at 8000 samples.}
    \label{fig:GOrderConvergence}
\end{figure}

Notice in Figures \ref{fig:Gc1Convergence} through \ref{fig:Gc3Convergence} that the graphs of ${\cal G}{12}$, ${\cal G}{13}$, ${\cal G}{14}$, and ${\cal G}{15}$ are practically identical, providing evidence for the validity of (\ref{eq:Gn-convergence-6}). To better illustrate this, sensitivity indices for all parameters and orders were retrieved at 8000 samples from Figures \ref{fig:Gc1Convergence}, \ref{fig:Gc2Convergence}, and \ref{fig:Gc3Convergence} and replotted, this time with the ${\cal G}$ orders on the horizontal axis (see Figure \ref{fig:GOrderConvergence}). Evidently, for orders greater than 11, the sensitivity indices for $c_1, c_2,$ and $c_3$ show negligible variation.

\begin{figure}[H]
    \centering
    \includegraphics{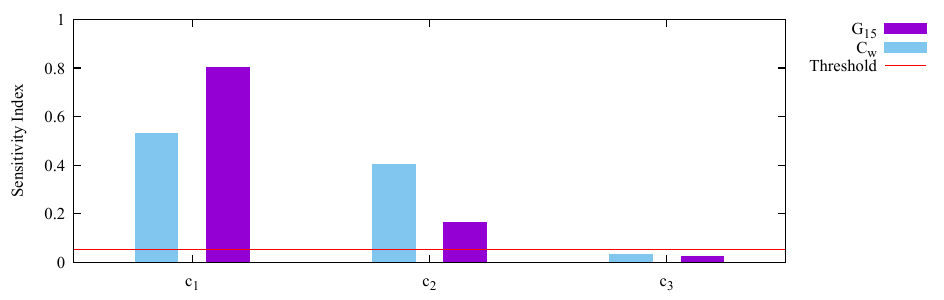}
    \caption{${\cal G}_{15}$-based sensitivity indices comparison to $C_w$-based sensitivity indices at 8000 samples}
    \label{fig:Gresults}
\end{figure}

Further insights were gained by retrieving the sensitivity indices of the three parameters at order 15 and plotting them in comparison to their $C_w$-based counterparts (see Figure \ref{fig:Gresults}). The results indicate agreement between ${\cal G}$ and $C_w$ with respect to the threshold. Additionally, the relative importance of $c_1, c_2,$ and $c_3$ is maintained during the transition between $C_w$ and ${\cal G}$. In both cases, it holds true that $SI_{c_1}>SI_{c_2}>SI_{c_3}.$

\subsection{Summary of SA Results}\label{sec:summaryOfSAresults}

\par
\begin{figure}[H]
    \centering
    \includegraphics{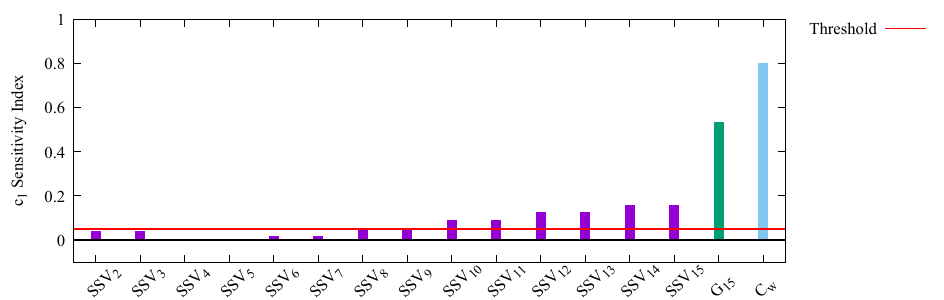}
    \caption{Sensitivity indices of $c_1$ parameter based on SSV of orders varying from 2 to 15, ${\cal G}$ of order 15, and $C_w$}
    \label{fig:SSV_G_Cw_c1}
\end{figure}

\begin{figure}[H]
    \centering
    \includegraphics{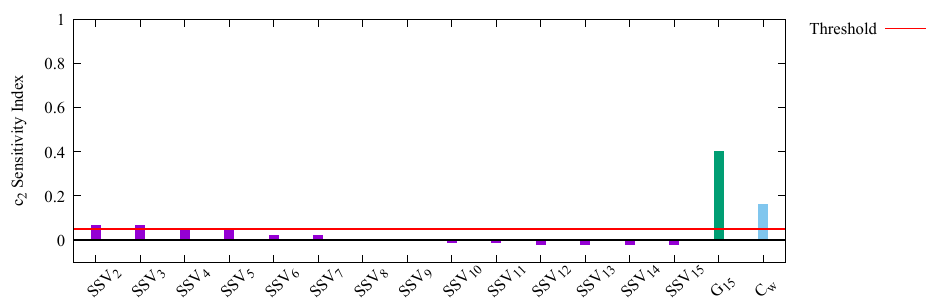}
    \caption{Sensitivity indices of $c_2$ parameter based on SSV of orders varying from 2 to 15, ${\cal G}$ of order 15, and $C_w$}
    \label{fig:SSV_G_Cw_c2}
\end{figure}

\begin{figure}[H]
    \centering
    \includegraphics{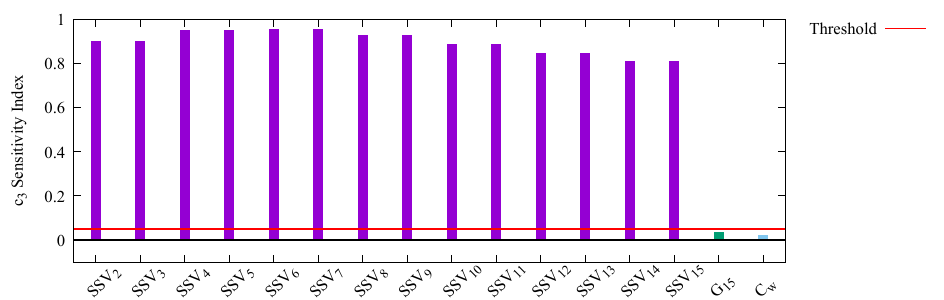}
    \caption{Sensitivity indices of $c_3$ parameter based on SSV of orders varying from 2 to 15, ${\cal G}$ of order 15, and $C_w$.}
    \label{fig:SSV_G_Cw_c3}
\end{figure}

\noindent The resulting sensitivity indices with respect to SSV, ${\cal G}$, and $C_w$ for parameters $c_1$, $c_2$, and $c_3$ are summarized in Figures \ref{fig:SSV_G_Cw_c1}, \ref{fig:SSV_G_Cw_c2}, and \ref{fig:SSV_G_Cw_c3}, respectively. Notice that in all figures, the behaviour of ${\cal G}_{15}$ is significantly more aligned with $C_w$ compared to that of SSV$_n$, especially in the case of $c_3$.

\par
\begin{figure}[H]
    \centering
    \includegraphics{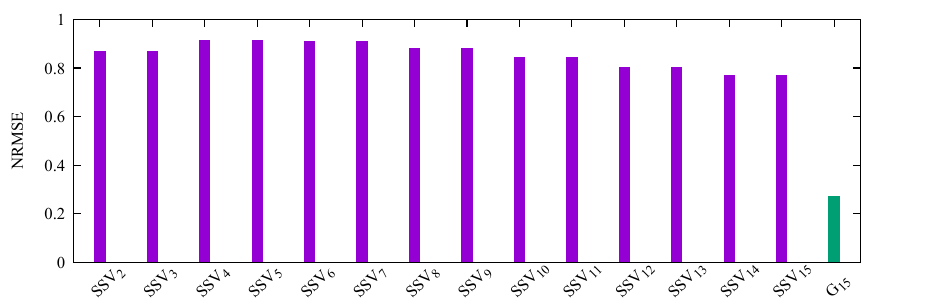}
    \caption{Error (\ref{eq:NRMSE}) between SSV- and $C_w$-based sensitivity indices compared to error between ${\cal G}$-based and $C_w$-based sensitivity indices.}
    \label{fig:NMRSE}
\end{figure}

\begin{figure}[H]
    \centering
    \includegraphics{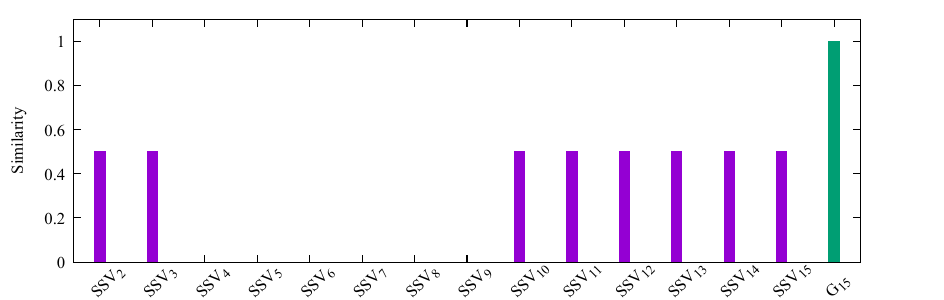}
    \caption{Similarity (\ref{eq:similarity}) between SSV- and $C_w$-based sensitivity indices compared to similarity between ${\cal G}$- and $C_w$-based sensitivity indices.}
    \label{fig:similarity}
\end{figure}

\noindent The correlation measures (\ref{eq:NRMSE}) and (\ref{eq:similarity}) can be applied to Figures \ref{fig:SSV_G_Cw_c1} through \ref{fig:SSV_G_Cw_c3} so that the performance of SSV relative to that of ${\cal G}$ can be more easily compared. Doing so for (\ref{eq:NRMSE}) and (\ref{eq:similarity}) results in Figures \ref{fig:NMRSE} and \ref{fig:similarity}, respectively. \textit{Notice in Figure \ref{fig:similarity} that the slender body operator ${\cal G}_{15}$ is able to capture 100\% of the sensitive parameters.}

\par\noindent Next, in Table \ref{tbl:finalsummary}, one can identify the sensitive parameters from Figures \ref{fig:SSV_G_Cw_c1}, \ref{fig:SSV_G_Cw_c2}, \ref{fig:SSV_G_Cw_c3}, in order of importance, as well as the normalised error and percentage similarity from Figures \ref{fig:NMRSE} and \ref{fig:similarity}. Notice that not only does ${\cal G}_{15}$ outperform SSV$_n$ in terms of error (\ref{eq:NRMSE}) and similarity (\ref{eq:similarity}), but also the relative importance of the parameters is preserved, so that $c_1$ can be correctly identified as more influential than $c_2$. These results illustrate that the inclusion of additional information regarding the underlying geometry with the ${\rm SSV}_n$ operators, in contrast to the restriction to moments of the form $M(N,0,0)$, as is the case with ${\cal G}$, produces results consistently less aligned with that of the physics operator $C_w$.\\

\begin{table}[h]
    \centering
    \begin{tabular}{ |c||c|c|c|  }
        \hline
        Quantity of Interest & Sensitive Parameters & NRMSE & Similarity \\
        \hline
        $C_w$ & $c_1 > c_2$ & N/A & N/A \\
        ${\cal G}_{15}$ & $c_1 > c_2$ & 0.267455 & 100\%  \\
        SSV$_{15}$ & $c_3 > c_1$ & 0.766259 & 50\% \\
        SSV$_{14}$ & $c_3 > c_1$ & 0.766271 & 50\% \\
        SSV$_{13}$ & $c_3 > c_1$ & 0.801278 & 50\% \\
        SSV$_{12}$ & $c_3 > c_1$ & 0.801279 & 50\% \\
        SSV$_{11}$ & $c_3 > c_1$ & 0.839496 & 50\% \\
        SSV$_{10}$ & $c_3 > c_1$ & 0.839486 & 50\% \\
        SSV$_{9}$ & $c_3$ & 0.877892 & 0\% \\
        SSV$_{8}$ & $c_3$ & 0.877874 & 0\% \\
        SSV$_{7}$ & $c_3$ & 0.908029 & 0\% \\
        SSV$_{6}$ & $c_3$ & 0.908006 & 0\% \\
        SSV$_{5}$ & $c_3$ & 0.911654 & 0\% \\
        SSV$_{4}$ & $c_3$ & 0.911578 & 0\% \\
        SSV$_{3}$ & $c_3 > c_2$ & 0.865169 & 50\% \\
        SSV$_{2}$ & $c_3 > c_2$ & 0.864787 & 50\% \\
        \hline
    \end{tabular}
    \caption{Comparison among $C_W$, SSV, and ${\cal G}$-based sensitivity analyses. The second column displays the order of sensitive parameters concerning the respective quantity of interest. The third and fourth columns present the error (\ref{eq:NRMSE}) and discrete similarity (\ref{eq:similarity}) between sensitivity indices evaluated with respect to ${\cal G}_{15}$ and $\text{SSV}_n$ compared to the indices obtained with $C_w$.}
    \label{tbl:finalsummary}
\end{table}
\subsection{Computational cost}
\begin{figure}[H]
    \centering
    \includegraphics{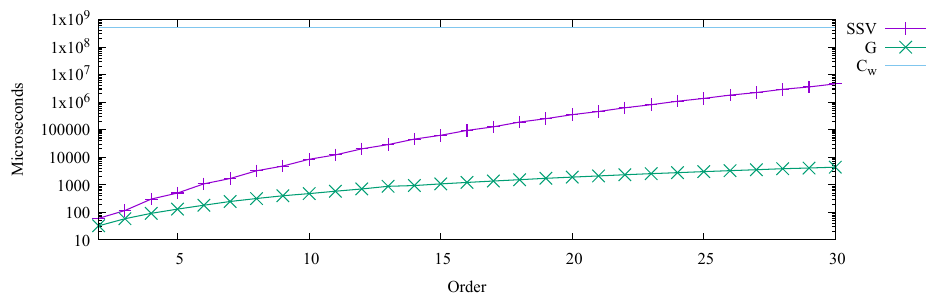}
    \caption{Per-run computational cost of evaluating $C_w$, SSV and ${\cal G}$ verses the order of SSV and ${\cal G}$.}
    \label{fig:GSSVtime}
\end{figure}

The computational cost to perform sensitivity analysis with geometric operators is significantly less than performing these analysis with $C_w$. On a PC with Intel(R) Xeon(R) Gold 6226 CPU with 2.70 GHz and 2.69 GHz processors and 128 GB of memory on average, each run of the IGA-BEM solver to evaluate $C_w$ took about 540 seconds, totalling around 210 hours to produce Figure \ref{fig:CwConvergence}. In contrast, both the SSV-based analysis and the ${\cal G}$-based analysis required computational time in the order of seconds. To compare the per-run computational cost of $C_w$ versus SSV (\ref{eq:SSV}) and ${\cal G}$ (\ref{eq:Gn}), a benchmark across 50 random samples was performed for all operators to establish, for a given order, the necessary computational cost. These results can be seen in Figure \ref{fig:GSSVtime}, where one can readily notice that for the tested orders, SSV is at least four orders of magnitude faster than $C_w$, and ${\cal G}$ is about 1 to 2 orders of magnitude faster than SSV. It is noteworthy that as the order increases, the computational cost of SSV rises much faster than that of ${\cal G}$, which can be attributed to the number of new moments that each order introduces.

\section{Conclusion and Future Work}\label{sec:6}

\par\noindent This study proposes a geometry-based operator to support computationally demanding physical models for reducing the dimensionality of the problem through sensitivity analysis. A general framework for matching physics-based quantities ${\cal P}$ to geometry-based ones is outlined in Figure \ref{fig:MethodologyDiagram}. Two geometry-based operators are proposed: the SSV, previously tested in \cite{GMGSA}, and the slender body operator ${\cal G}$, derived in this paper based on slender body theory. Choosing ${\cal P} = C_w$ as the wave resistance coefficient and ${\cal D}$ as the modified Wigley hull parametric modeller, the framework is applied twice: once to investigate the compatibility between $C_w$ and SSV, and once for $C_w$ and ${\cal G}$.

\par\noindent The results indicate that ${\cal G}$ outperforms SSV significantly, as summarised in Table \ref{tbl:finalsummary}. In this table, one can observe that the slender body operator outperforms SSV for all displayed orders, achieving 100\% similarity to $C_w$ while also preserving the relative importance of the sensitive parameters. Furthermore, the computational cost of the ${\cal G}$-based operator was significantly less than that of $C_w$, with the former taking seconds and the latter taking days. Additionally, as shown in Figure \ref{fig:GSSVtime}, the computational cost of ${\cal G}$ is significantly lower than that of SSV, a difference that becomes more pronounced as the order increases.

\par\noindent While the slender body operator has proven to be a promising candidate geometric operator for the problem of wave resistance when paired with the parametric modeller (\ref{eq:wigley}), several avenues have not been explored in this study. First and foremost, the framework of Figure \ref{fig:MethodologyDiagram} was applied only for one parametric modeller, which also has a small number of parameters. Further experimentation with richer parametric modellers is necessary to expand the applicability of the pair ${\cal G}-{\cal P}$. The two parametric modellers tested in \cite{GMGSA}, for which SSV performed satisfactorily, are promising candidates in this regard, allowing for a more in-depth investigation into the seemingly better performance of $\cal G$ compared to SSV.

Next, the relatively short computational times of the ${\cal G}$-based analysis can be attributed in part to the closed-form expression for the moments (\ref{eq:MTpqr-final-FULL}), made possible by the simple analytical expression of (\ref{eq:wigley}). However, it is acknowledged that this may not be possible with more complicated geometries. Nevertheless, relevant experimentation has shown that the increase in computational time is not drastic in the general-surface case, still providing a significant reduction when compared to the $C_w$-based sensitivity analysis.

Finally, there is potential to enhance ${\cal G}$ further by considering the generalised formulation of Vossers' integral, as provided in \cite{maruo1962calculation}. This has the potential to increase the performance of the slender body operator.

\section*{ACKNOWLEDGEMENTS}
This work received funding from the European Union's Horizon 2020 research and innovation programme under the Marie Skłodowska-Curie grant agreement No 860843, PI for the University of Strathclyde: P.D. Kaklis.

\bibliographystyle{elsarticle-num}
\bibliography{ref}

\begin{thebibliography}{10}
\bibitem{maruo1962calculation}
H.~Maruo, Calculation of the wave resistance of ships, the draught of which is
  as small as the beam, Journal of Zosen Kiokai 1962~(112) (1962) 21--37.
\end{thebibliography}
\appendix 
\section{Convergence of the Slender Body Operator}\label{appendix:slenderbodyconvergentseries}

\par\noindent In this section we show the convergence of (\ref{eq:Gn}), for which a comparison and a ratio test suffice. For $S(x)$ the sectional area curve of the hull at the longitudinal position $x$ assume that there exists an $M>0$ with,
\begin{gather}
    |S''(x)| \le M,\quad \forall x.
\end{gather}
Then,
\begin{align}\label{eq:Gn-convergence-1}
    S''(x)S''(\xi)(x-\xi)^{2k} \le M^2(x-\xi)^{2k},\quad \forall x,\ \xi.
\end{align}
Looking at (\ref{eq:Ik-0}) and since the right-hand side of (\ref{eq:Gn-convergence-1}) is certainly integrable,
\begin{align}\label{eq:Gn-convergence-1p2}
    I_k &= \int_{-L/2}^{L/2}\int_{-L/2}^{L/2}S''(x)S''(\xi)(x-\xi)^{2k}dxd\xi \nonumber\\
    &\le M^2\int_{-L/2}^{L/2}\int_{-L/2}^{L/2}(x-\xi)^{2k}dxd\xi\nonumber\\
    &=M^2\int_{-L/2}^{L/2}\Bigg[\frac{(x-\xi)^{2k+1}}{2k+1}\Bigg]_{-L/2}^{L/2}d\xi\nonumber\\
    &=\frac{M^2}{2k+1}\Bigg[\frac{-(L/2-\xi)^{2k+2}+(-L/2-\xi)^{2k+2}}{2k+2}\Bigg]_{-L/2}^{L/2}\nonumber\\
    &=\frac{2M^2L^{2k+2}}{(2k+1)(2k+2)}.
\end{align}
Using (\ref{eq:Gn-convergence-1p2}) and looking at (\ref{eq:Gn}),
\begin{align}\label{eq:Gn-convergence-1p5}
    |f(k;K)I_k|\le|f(k;K)\frac{2M^2L^{2k+2}}{(2k+1)(2k+2)}|,
\end{align}
which means that it is enough to show that the following series is absolutely convergent,
\begin{align}
    \sum_{k=0}^{\infty}\beta_k, \quad \beta_k=f(k;K)\frac{2M^2L^{2k+2}}{(2k+1)(2k+2)}.
\end{align}
This can be done by application of the ratio test,
\begin{align}\label{eq:Gn-convergence-2}
    \frac{|\beta_{k+1}|}{|\beta_k|} &= \frac{f(k+1;K)}{f(k;K)}L^2\frac{(2k+1)(2k+2)}{(2k+3)(2k+4)}\nonumber\\
    &=\frac{|f(k+1;K)|}{|f(k;K)|}O(1),\ k\to \infty.
\end{align}
Then, recalling the definition of $f(k;K)$ from (\ref{eq:fkK}),
\begin{align}\label{eq:Gn-convergence-3}
    \frac{|f(k+1;K)|}{|f(k;K)|}&=\frac{K^22^{2k-1}(k!)^2}{2^{2k+1}((k+1)!)^2}\cdot\frac{ln(0.5K)+\gamma - h(k+1)}{ln(0.5K)+\gamma - h(k)}\nonumber\\
    &=\frac{K^2}{4(k+1)^2}\cdot\frac{ln(0.5K)+\gamma - h(k+1)}{ln(0.5K)+\gamma - h(k)}\nonumber\\
    &=\frac{ln(0.5K)+\gamma - h(k+1)}{ln(0.5K)+\gamma - h(k)}O(k^{-2}),\ k\to \infty.
\end{align}
Now, for the harmonic series $h(k)$, it is true that $h(k)=O(\ln(k)),\ k\to \infty$. In light of this,
\begin{align}\label{eq:Gn-convergence-4}
    \frac{\ln(0.5K)+\gamma - h(k+1)}{\ln(0.5K)+\gamma - h(k)} &= \frac{h(k+1)}{h(k)}\cdot\frac{\displaystyle1-\frac{\ln(0.5K)+\gamma}{h(k+1)}}{\displaystyle1-\frac{\ln(0.5K)+\gamma}{h(k)}}\nonumber\\
    &=\frac{\displaystyle\sum_{i=0}^{k+1}\frac{1}{i}}{\displaystyle\sum_{i=0}^k\frac{1}{i}}\cdot\frac{1-O((\ln(k+1))^{-1})}{1-O((\ln(k))^{-1})}\nonumber\\
    &=\Bigg(1+\frac{1}{(k+1)h(k)}\Bigg)\cdot\frac{1+o(1)}{1+o(1)}\nonumber\\
    &=(1+O(k^{-1}(\ln(k))^{-1}))(1+o(1))\nonumber\\
    &=1+o(1),\ k\to \infty.
\end{align}
Substituting (\ref{eq:Gn-convergence-4}) into (\ref{eq:Gn-convergence-3}),
\begin{align}\label{eq:Gn-convergence-5}
    \frac{|f(k+1;K)|}{|f(k; K)|}=(1+o(1))O(k^{-2})=o(1),\ k\to\infty.
\end{align}
Substituting (\ref{eq:Gn-convergence-5}) into (\ref{eq:Gn-convergence-2}),
\begin{align}
    \frac{|\beta_{k+1}|}{|\beta_k|} &= o(1),\ k\to\infty,
\end{align}
which means that $\{\beta_i\}$ is an absolutely convergent series. By comparison, (\ref{eq:Gn-convergence-1p5}) implies that
\begin{gather}\label{eq:Gn-convergence-6}
    \lim_{n\to\infty}{\cal G}(n)\text{ exists and}\in\mathbb{R}.
\end{gather}

\section{Closed Forms of the Geometric Moments of the Modified Wigley Hull}\label{appendix:closedformsappendix}
\par\noindent In this section, we express the moments (\ref{eq:Mpqr}) and (\ref{eq:MTpqr}) for (\ref{eq:wigley}) analytically. We begin with $M(p,q,r)$ by substituting (\ref{eq:wigley}) into (\ref{eq:Mpqr}),
\begin{align}\label{eq:Mpqr-1}
    M(p,q,r) &= \int_{-L/2}^{L/2}\int_{0}^T\int_{-B\eta/2}^{B\eta/2}x^py^qz^rdydzdx\nonumber\\
    &= \int_{-L/2}^{L/2}\int_{0}^Tx^pz^r\Bigg[\frac{y^{q+1}}{q+1}\Bigg]_{-B\eta/2}^{B\eta/2}dzdx\nonumber\\
    &= \Big(\frac{B}{2}\Big)^{q+1}\frac{(1-(-1)^{q+1})}{q+1}\int_{-L/2}^{L/2}\int_{0}^Tx^pz^r\eta^{q+1}dzdx.
\end{align}
Notice that if $q$ is odd then $M(p,q,r)=0$ which can be attributed to the breadth-wise symmetry of the Wigley hull. Assuming $q$ to be even, we write
\begin{align}\label{eq:Mpqr-2}
    M(p,q,r) = \Big(\frac{B}{2}\Big)^{q+1}\frac{2}{q+1}\int_{-L/2}^{L/2}\int_{0}^Tx^pz^r\eta^{q+1}dzdx,\quad q = \text{even}.
\end{align}
To proceed, we expand $\eta^{q+1}$,
\begin{align}
    \eta^{q+1} &= \Big((1-\zeta^2)(1-\xi^2)(1+c_1\xi^2+c_2\xi^4)+c_3\zeta^2(1-\zeta^8)(1-\xi^2)^4\Big)^{q+1}\nonumber\\
    &=\sum_{i_1=0}^{q+1}\binom{q+1}{i_1}\Big((1-\zeta^2)(1-\xi^2)(1+c_1\xi^2+c_2\xi^4)\Big)^{q+1-i_1}\Big(c_3\zeta^2(1-\zeta^8)(1-\xi^2)^4\Big)^{i_1}\nonumber\\
    &=\sum_{i_1=0}^{q+1}\binom{q+1}{i_1}(1-\zeta^2)^{q+1-i_1}(1-\xi^2)^{q+1+3i_1}(1+c_1\xi^2+c_2\xi^4)^{q+1-i_1}c_3^{i_1}\zeta^{2i_1}(1-\zeta^8)^{i_1}\nonumber\\
    &=\sum_{i_1=0}^{q+1}\binom{q+1}{i_1}\Bigg(\sum_{i_2=0}^{q+1-i_1}\binom{q+1-i_1}{i_2}(-1)^{i_2}\zeta^{2i_2}\Bigg)\Bigg(\sum_{i_3=0}^{q+1+3i_1}\binom{q+1+3i_1}{i_3}(-1)^{i_3}\xi^{2i_3}\Bigg)\nonumber\\
    &\quad\cdot\Bigg(\sum_{i_4=0}^{q+1-i_1}\binom{q+1-i_1}{i_4}(c_1\xi^2+c_2\xi^4)^{i_4}\Bigg)c_3^{i_1}\zeta^{2i_1}\Bigg(\sum_{i_6=0}^{i_1}\binom{i_1}{i_6}(-1)^{i_6}\zeta^{8i_6}\Bigg)\nonumber\\
    &=\sum_{i_1=0}^{q+1}\binom{q+1}{i_1}\Bigg(\sum_{i_2=0}^{q+1-i_1}\binom{q+1-i_1}{i_2}(-1)^{i_2}\zeta^{2i_2}\Bigg)\Bigg(\sum_{i_3=0}^{q+1+3i_1}\binom{q+1+3i_1}{i_3}(-1)^{i_3}\xi^{2i_3}\Bigg)\nonumber\\
    &\quad\cdot\Bigg(\sum_{i_4=0}^{q+1-i_1}\binom{q+1-i_1}{i_4}\sum_{i_5=0}^{i_4}\binom{i_4}{i_5}c_1^{i_5}\xi^{2i_5}c_2^{i_4-i_5}\xi^{4i_4-4i_5}\Bigg)c_3^{i_1}\zeta^{2i_1}\Bigg(\sum_{i_6=0}^{i_1}\binom{i_1}{i_6}(-1)^{i_6}\zeta^{8i_6}\Bigg).
\end{align}
Grouping together all the sum-operators and combining common factors,
\begin{align}\label{eq:etap-final}
    \eta^{q+1} &= \sum_{i_1=0}^{q+1}\ \sum_{i_2=0}^{q+1-i_1}\ \sum_{i_3=0}^{q+1+3i_1}\ \sum_{i_4=0}^{q+1-i_1}\ \sum_{i_5=0}^{i_4}\ \sum_{i_6=0}^{i_1}\ \binom{q+1}{i_1}\binom{q+1-i_1}{i_2}\binom{q+1+3i_1}{i_3}\nonumber\\&\quad\cdot\binom{q+1-i_1}{i_4}\binom{i_4}{i_5}\binom{i_1}{i_6}(-1)^{i_2+i_3+i_6}c_1^{i_5}c_2^{i_4-i_5}c_3^{i_1}\zeta^{2i_1+2i_2+8i_6}\xi^{2i_3+4i_4-2i_5}.
\end{align}
Substituting (\ref{eq:etap-final}) into (\ref{eq:Mpqr-2}) and factoring all constants out of the double integral,
\begin{align}\label{eq:Mpqr-3}
    M(p,q,r) &= \Big(\frac{B}{2}\Big)^{q+1}\frac{2}{q+1} \sum_{i_1=0}^{q+1}\ \sum_{i_2=0}^{q+1-i_1}\ \sum_{i_3=0}^{q+1+3i_1}\ \sum_{i_4=0}^{q+1-i_1}\ \sum_{i_5=0}^{i_4}\ \sum_{i_6=0}^{i_1}\Bigg[\nonumber\\&\quad \binom{q+1}{i_1}\binom{q+1-i_1}{i_2}\binom{q+1+3i_1}{i_3}\binom{q+1-i_1}{i_4}\binom{i_4}{i_5}\binom{i_1}{i_6}\nonumber\\
&\quad\cdot(-1)^{i_2+i_3+i_6}c_1^{i_5}c_2^{i_4-i_5}c_3^{i_1}\int_{-L/2}^{L/2}\xi^{2i_3+4i_4-2i_5}x^pdx\int_{0}^T\zeta^{2i_1+2i_2+8i_6}z^rdz\Bigg],\quad q = \text{even}.
\end{align}
Since $\xi = 2x/L$,
\begin{align}
    \int_{-L/2}^{L/2}\xi^{2i_3+4i_4-2i_5}x^pdx &=(2/L)^{2i_3+4i_4-2i_5}\int_{-L/2}^{L/2}x^{p+2i_3+4i_4-2i_5}dx\nonumber\\
    &=(2/L)^{2i_3+4i_4-2i_5}\Bigg[\frac{x^{p+2i_3+4i_4-2i_5+1}}{p+2i_3+4i_4-2i_5+1}\Bigg]_{-L/2}^{L/2}\nonumber\\
    &=\Big(\frac{L}{2}\Big)^{p+1}\frac{(1-(-1)^{p+1})}{p+2i_3+4i_4-2i_5+1}.
\end{align}
Again, if $p$ is odd, $M(p,q,r)=0$ which can be attributed to the length-wise symmetry of the Wigley hull. Assuming that $p$ is even,
\begin{align}\label{eq:xiIntegral}
    \int_{-L/2}^{L/2}\xi^{2i_3+4i_4-2i_5}x^pdx = \Big(\frac{L}{2}\Big)^{p+1}\frac{2}{p+2i_3+4i_4-2i_5+1},\quad p=\text{even}.
\end{align}
Similarly for $\zeta = z/T$,
\begin{align}\label{eq:zetaIntegral}
    \int_{0}^T\zeta^{2i_1+2i_2+8i_6}z^rdz = T^{r+1}\frac{1}{r+2i_1+2i_2+8i_6+1}.
\end{align}
Substituting (\ref{eq:xiIntegral}) and (\ref{eq:zetaIntegral}) into (\ref{eq:Mpqr-3}),
\begin{align}\label{eq:Mpqr-4}
     M(p,q,r) &= \Big(\frac{L}{2}\Big)^{p+1}\Big(\frac{B}{2}\Big)^{q+1}T^{r+1}\frac{4}{q+1} \sum_{i_1=0}^{q+1}\ \sum_{i_2=0}^{q+1-i_1}\ \sum_{i_3=0}^{q+1+3i_1}\ \sum_{i_4=0}^{q+1-i_1}\ \sum_{i_5=0}^{i_4}\ \sum_{i_6=0}^{i_1}\Bigg[\nonumber\\&\quad \binom{q+1}{i_1}\binom{q+1-i_1}{i_2}\binom{q+1+3i_1}{i_3}\binom{q+1-i_1}{i_4}\binom{i_4}{i_5}\binom{i_1}{i_6}\nonumber\\
&\quad\cdot(-1)^{i_2+i_3+i_6}c_1^{i_5}c_2^{i_4-i_5}c_3^{i_1}(p+2i_3+4i_4-2i_5+1)^{-1}(r+2i_1+2i_2+8i_6+1)^{-1}\Bigg],\\
&\nonumber\text{where } p,\ q = \text{even}.
\end{align}
The presentation of (\ref{eq:Mpqr-4}) can be simplified considerably by noticing that the sums can be split into three groups: the first group consists only of the sum with index $i_1$, the second group of the sums with indexes $\{i_3,\ i_4,\ i_5\}$ and the final group with $\{i_2,\ i_6\}$. Then,
\begin{subequations}\label{eq:Mpqr-final-FULL}
\begin{align}\label{eq:Mpqr-final}
    M(p,q,r) = \begin{cases}\displaystyle
        \Big(\frac{L}{2}\Big)^{p+1}\Big(\frac{B}{2}\Big)^{q+1}T^{r+1}\frac{4}{q+1} \sum_{i_1=0}^{q+1}\binom{q+1}{i_1}c_3^{i_1}{\cal X}_{i_1}{\cal Z}_{i_1},\quad \text{$p,\ q$ are even}\\
        0,\quad \text{else},
    \end{cases}
\end{align}
where
\begin{align}\label{eq:Xi}
    {\cal X}_{i_1} = \sum_{i_3=0}^{q+1+3i_1}\ \sum_{i_4=0}^{q+1-i_1}\ \sum_{i_5=0}^{i_4}\binom{q+1+3i_1}{i_3}\binom{q+1-i_1}{i_4}\binom{i_4}{i_5}(-1)^{i_3}c_1^{i_5}c_2^{i_4-i_5}(p+2i_3+4i_4-2i_5+1)^{-1},
\end{align}
and 
\begin{align}\label{eq:Zi-standard}
    {\cal Z}_{i_1} = \sum_{i_2=0}^{q+1-i_1}\ \sum_{i_6=0}^{i_1}\binom{q+1-i_1}{i_2}\binom{i_1}{i_6}(-1)^{i_2+i_6}(r+2i_1+2i_2+8i_6+1)^{-1}.
\end{align}
\end{subequations}
Having (\ref{eq:Mpqr-final}) it is easy to evaluate an analogous expression for the translation invariant moment $M_T(p,q,r)$. as given in (\ref{eq:MTpqr}),
\begin{gather}\label{eq:MTpqr-1}
    M_T(p,q,r) = \int_{-L/2}^{L/2}\int_{0}^T\int_{-B\eta/2}^{B\eta/2}(x-C_x)^p(y-C_y)^q(z-C_z)^rdydzdx.
\end{gather}
First, notice that by (\ref{eq:Mpqr-final}), $M(1,0,0)=M(0,1,0)=0$. Then, looking at (\ref{eq:centroid}), $C_x=C_y=0$,
\begin{gather}\label{eq:MTpqr-2}
    M_T(p,q,r) = \int_{-L/2}^{L/2}\int_{0}^T\int_{-B\eta/2}^{B\eta/2}x^py^q(z-C_z)^rdydzdx.
\end{gather}
In the derivation of (\ref{eq:Mpqr-final}), the $z^r$ factor (now $(z-C_z)^r$) was integrated at (\ref{eq:zetaIntegral}). One can then easily be convinced that the only difference between $M_T(p,q,r)$ and (\ref{eq:Mpqr-final}) will be in ${\cal Z}_{i_1}$. Specifically, looking at (\ref{eq:zetaIntegral})
\begin{align}\label{eq:zetaIntegral-Tinvariance}
    \int_{0}^T\zeta^{2i_1+2i_2+8i_6}(z-C_z)^rdz &= \Big(T^{2i_1+2i_2+8i_6}\Big)^{-1}\int_0^{T}z^{2i_1+2i_2+8i_6}\sum_{j = 0}^r\binom{r}{j}z^{r-j}(-1)^jC_z^jdz\nonumber\\
    &=\Big(T^{2i_1+2i_2+8i_6}\Big)^{-1}\sum_{j = 0}^r\binom{r}{j}(-1)^jC_z^j\Bigg[\frac{z^{r+2i_1+2i_2+8i_6-j+1}}{r+2i_1+2i_2+8i_6-j+1}\Bigg]_0^T\nonumber\\
    &=T^{r+1}\sum_{j = 0}^r\binom{r}{j}(-1)^j\Big(\frac{C_z}{T}\Big)^j(r+2i_1+2i_2+8i_6-j+1)^{-1}.
\end{align}
So that finally,
\begin{subequations}\label{eq:MTpqr-final-FULL}
\begin{align}\label{eq:MTpqr-final}
    M_T(p,q,r) = \begin{cases}\displaystyle
        \Big(\frac{L}{2}\Big)^{p+1}\Big(\frac{B}{2}\Big)^{q+1}T^{r+1}\frac{4}{q+1} \sum_{i_1=0}^{q+1}\binom{q+1}{i_1}c_3^{i_1}{\cal X}_{i_1}{\cal Z}_{i_1},\quad \text{$p,\ q$ are even}\\
        0,\quad \text{else},
    \end{cases}
\end{align}
with
\begin{align}\label{eq:Zi-transformationinvariant}
    {\cal Z}_{i_1} = \sum_{i_2=0}^{q+1-i_1}\ \sum_{i_6=0}^{i_1}\ \sum_{j = 0}^r\binom{q+1-i_1}{i_2}\binom{i_1}{i_6}\binom{r}{j}(-1)^{i_2+i_6+j}\Big(\frac{C_z}{T}\Big)^j(r+2i_1+2i_2+8i_6-j+1)^{-1},
\end{align}
\end{subequations}
where $C_z$ is given by (\ref{eq:centroid}) and ${\cal X}_{i_1}$ by (\ref{eq:Xi}). It should be noted that the validity of equations (\ref{eq:Mpqr-final-FULL}) and (\ref{eq:MTpqr-final-FULL}) were verified via the computer algebraic system in MAPLE$^{\copyright}$ \footnote{https://www.maplesoft.com/products/Maple/}.

\section{Closed Form of the Sectional Area Curve of the Modified Wigley Hull}\label{appendix:wigley_sac}

In this section, a closed form for the sectional area of (\ref{eq:wigley}) at a given longitudinal position $x\in[-L/2,L/2]$, $S(x)$, is derived. We can directly write
\begin{align}\label{eq:sac_1}
    S(x) =\int_0^T\int_{-B\eta/2}^{B\eta/2}dydz=B\int_0^T\eta dz = B(1-\xi^2)(1+c_1\xi^2+c_2\xi^4)\int_0^T(1-\zeta^2)dz+Bc_3(1-\xi^2)^4\int_0^T(\zeta^2-\zeta^{10})dz
\end{align}
where $\xi = 2x/L$ and $\zeta = z/T$. Evaluating the integrals in (\ref{eq:sac_1}),
\begin{align}\label{eq:sac}
    S(x) = BT\Bigg(\frac{2}{3}(1-\xi^2)(1+c_1\xi^2+c_2\xi^4)+\frac{8c_3}{33}(1-\xi^2)^4\Bigg)
\end{align}
Differentiating (\ref{eq:sac}) with respect to $x$,
\begin{align}\label{eq:dsac}
    S'(x) = \frac{2BT}{L}\Bigg(-\frac{4\xi}{3}(1+c_1\xi^2+c_2\xi^4)+\frac{2}{3}(1-\xi^2)(2c_1\xi+4c_2\xi^3-32 c_3\xi(1-\xi^2)^2/11)\Bigg)
\end{align}
For $x=\pm L/2$, $\xi(x=\pm L/2)=\pm 1$ we get,
\begin{align}\label{eq:dsac_1}
    S'(\pm L/2) = \mp \frac{8BT}{3L}(1+c_1+c_2)
\end{align}
which attains its maximum at $c_1=c_2=1$, $|S'(\pm L/2;c_1=c_2=1)| = 8BT/L$
\end{document}